\documentclass[10pt]{article}
\usepackage{amssymb,amsfonts}
\usepackage{amsmath}
\usepackage{amscd,amssymb,amsthm}
\usepackage{graphics}
\usepackage{graphicx}
\usepackage{hyperref}

\graphicspath{ {./images/} }
\hypersetup{
    colorlinks=true,citecolor=blue
  }
\newtheorem{theorem}{Theorem}[section]

\newtheorem{corollary}[theorem]{Corollary}

\newtheorem{definition}[theorem]{Definition}

\newtheorem{proposition}[theorem]{Proposition}
\newtheorem{remark}[theorem]{Remark}

\begin{document}

\title{Analysis of a class of Kolmogorov systems}
\author{G. Moza, C. Lazureanu\thanks{
Department of Mathematics, Politehnica University of Timisoara, Romania;
gheorghe.tigan@upt.ro}, F. Munteanu, C. Sterbeti, A. Florea\thanks{%
Department of Applied Mathematics, University of Craiova, Romania} }
\date{}
\maketitle

\begin{abstract}
A two-dimensional Kolmogorov system with two parameters
and having a degenerate condition is studied in this work. We obtain local
analytical properties of the system when the parameters vary in a
sufficiently small neighborhood of the origin. The behavior of the system is
described by bifurcation diagrams. Applications of Kolmogorov systems can be
found particularly in modeling population dynamics in biology and ecology.
\end{abstract}

%cristian.lazureanu@upt.ro
%munteanufm@gmail.com
%aurelia\_florea@yahoo.com
%sterbetiro@yahoo.com

\section{Introduction}

In \cite{GT1} we have studied a class of two-dimensional Kolmogorov systems 
\cite{jl1} of the form 
\begin{equation}
\left\{ 
\begin{tabular}{lll}
$\frac{d\xi _{1}}{dt}$ & $=$ & $\xi _{1}\left( \mu _{1}-\theta \left( \mu
\right) \xi _{1}+\gamma \left( \mu \right) \xi _{2}-M\left( \mu \right) \xi
_{1}\xi _{2}+N\left( \mu \right) \xi _{1}^{2}\right) $ \\ 
&  &  \\ 
$\frac{d\xi _{2}}{dt}$ & $=$ & $\xi _{2}\left( \mu _{2}-\delta \left( \mu
\right) \xi _{1}+\xi _{2}+S\left( \mu \right) \xi _{1}^{2}+P\left( \mu
\right) \xi _{2}^{2}\right) $%
\end{tabular}%
\right.  \label{fn2}
\end{equation}%
in a non-degenerate context given by $\theta \left( 0\right) \delta \left(
0\right) \neq 0$ and $\gamma \left( 0\right) <0.$ The coefficients $\theta
\left( \mu \right) ,$ $\gamma \left( \mu \right) ,$ $\delta \left( \mu
\right) ,$ $M(\mu ),$ $N(\mu ),$ $S(\mu )$ and $P(\mu ),$ are smooth
functions of class $C^{k},$ $k\geq 1,$ and the parameter $\mu =\left( \mu
_{1},\mu _{2}\right) \in 
%TCIMACRO{\U{211d} }%
%BeginExpansion
\mathbb{R}
%EndExpansion
^{2}$ is chosen such that $\left\vert \mu \right\vert =\sqrt{\mu
_{1}^{2}+\mu _{2}^{2}}$ is sufficiently small, $0\leq \left\vert \mu
\right\vert \ll 1;$ for brevity we denote $\left( 0,0\right) $ by $0.$

We are concerned in this work with local properties of the system \eqref{fn2}
in a degenerate framework given by 
\begin{equation}
\theta \left( 0\right) \delta \left( 0\right) =0\text{ and }\gamma \left(
0\right) <0.  \label{dege}
\end{equation}%
This new condition changes significantly the behavior of the system. Since
many applications of Kolmogorov systems use positive variables \cite{adak}, 
\cite{green}, \cite{jl2}, the phase space of \eqref{fn2} which we consider
in this work is the first quadrant $\xi _{1}\geq 0$ and $\xi _{2}\geq 0.$
Thus, the equilibrium points will be studied only when their coordinates are
positive or zero. A similar study can be performed for the other quadrants.

Applications of Kolmogorov systems can be found mainly in biology \cite{adak}%
, \cite{bel}, \cite{rafi} and ecology \cite{kot}, \cite{may}, \cite{yang},
modeling population dynamics. A particular class of Kolmogorov systems is
the class of Lotka--Volterra systems \cite{bra}, \cite{cooke}, \cite{free}, 
\cite{lob}, \cite{precup}, \cite{yama}, which are widely used for modeling
the behavior of interacting biological different species of predator-prey
type. For example, an autonomous Lotka--Volterra type competitive system for 
$N-$species to control the gut microbiota by antibiotics has been presented
recently in \cite{dong}. Many environmental, engineering, economics and
mechanical models can be reduced to some kind of Kolmogorov systems \cite%
{yuan}.

\section{Analysis of the system when $\protect\theta(0) =0$ and $\protect%
\delta(0) \neq 0$}

Assume $\theta \left( 0\right) =0,$ $\delta \left( 0\right) \neq 0$ and $%
N\left( 0\right) \neq 0,$ where $\gamma \left( 0\right) <0.$ Then $\theta
\left( \mu \right) =\frac{\partial \theta \left( 0\right) }{\partial \mu _{1}%
}\mu _{1}+\frac{\partial \theta \left( 0\right) }{\partial \mu _{2}}\mu
_{2}+O\left( \left\vert \mu \right\vert ^{2}\right) $ and assume $\frac{%
\partial \theta \left( 0\right) }{\partial \mu _{2}}\overset{not}{=}\theta
_{2}\neq 0.$ Throughout this work, 
\begin{equation*}
O\left( \left\vert \mu \right\vert ^{k}\right) =\sum_{i+j\geq k}c_{ij}\mu
_{1}^{i}\mu _{2}^{j}
\end{equation*}%
denotes a Taylor series starting with terms of order at least $k\geq 1.$

\begin{remark}
To save symbols, we denote further by $\gamma \left( 0\right) =\gamma,$ $%
\delta \left( 0\right) =\delta ,$ $N\left( 0\right) =N,$ $M\left( 0\right)
=M,$ $P\left( 0\right) =P $ and $S\left( 0\right) =S.$
\end{remark}

Different to the case $\theta \left( 0\right) \delta \left( 0\right) \neq 0,$
now it may exist two equilibrium points $E_{11}\left( \xi _{11}\left( \mu
\right) ,0\right) $ and $E_{12}\left( \xi _{12}\left( \mu \right) ,0\right) $
lying on the $\xi _{1}-$axis, where

\begin{equation}
\xi _{11}\left( \mu \right) =\frac{1}{2N\left( \mu \right) }\left( \theta
\left( \mu \right) +\sqrt{\Delta \left( \mu \right) }\right) \text{ and }\xi
_{12}\left( \mu \right) =\frac{1}{2N\left( \mu \right) }\left( \theta \left(
\mu \right) -\sqrt{\Delta \left( \mu \right) }\right) ,  \label{xi12g}
\end{equation}%
provided that $\Delta \left( \mu \right) =\theta ^{2}\left( \mu \right)
-4N\left( \mu \right) \mu _{1}\geq 0.$ In its lowest terms, we can write%
\begin{equation}
\Delta \left( \mu \right) =\theta _{2}^{2}\mu _{2}^{2}\left( 1+O\left( \mu
_{2}\right) \right) -4N\mu _{1}\left( 1+O\left( \left\vert \mu \right\vert
\right) \right) .  \label{demu}
\end{equation}%
Since $N=N\left( 0\right) \neq 0,$ the bifurcation curve $\Delta \left( \mu
\right) =0$ exists and is unique in the parametric plane $\mu_1O\mu_2$ for
all $\mu $ with $0\leq \left\vert \mu \right\vert \ll 1.$ This result is
obtained from the Implicit Function Theorem (IFT) since $\Delta \left(
0,0\right) =0$ and $\frac{\partial \Delta }{\partial \mu _{1}}\left(
0,0\right) =-4N\neq 0.$ Moreover, it follows from \eqref{demu} and IFT that
the curve $\Delta \left( \mu \right) =0$ has the expression $4N\mu
_{1}=\theta _{2}^{2}\mu _{2}^{2}\left( 1+O\left( \mu _{2}\right) \right)$
for all $\left\vert \mu \right\vert$ sufficiently small.

Denoting by $\Delta $ the bifurcation curve $\Delta \left( \mu \right) =0,$
it follows that $\Delta $ in its lowest terms becomes 
\begin{equation}
\Delta =\left\{ \left( \mu _{1},\mu _{2}\right) \in 
%TCIMACRO{\U{211d} }%
%BeginExpansion
\mathbb{R}
%EndExpansion
^{2}\left\vert -4N\mu _{1}+\theta _{2}^{2}\mu _{2}^{2}=0\right. \right\} .
\label{de1}
\end{equation}

\begin{remark}
For the description of qualitative properties of solutions of the system %
\eqref{fn2}, $\Delta \left( \mu \right) $ from \eqref{demu} can be
approximated by $\Delta =-4N\mu _{1}+\theta _{2}^{2}\mu _{2}^{2},$
respectively, $\xi _{11}\left( \mu \right) $ and $\xi _{12}\left( \mu
\right) $ by 
\begin{equation}
\xi _{11}=\frac{1}{2N}\left( \theta _{2}\mu _{2}+\sqrt{\Delta }\right) \text{
and }\xi _{12}=\frac{1}{2N}\left( \theta _{2}\mu _{2}-\sqrt{\Delta }\right) .
\label{xi12}
\end{equation}
\end{remark}

Denote further by $\Delta _{+}=\Delta \cap \left\{ \mu _{2}>0\right\} $ and $%
\Delta _{-}=\Delta \cap \left\{ \mu _{2}<0\right\} $ the two branches of the
curve $\Delta $ lying in the upper, respectively, lower half-plane.

The two equilibria $E_{11}\left( \xi _{11},0\right) $ and $E_{12}\left( \xi
_{12},0\right) $ exist whenever $\Delta \geq 0.$ However, since the phase
space for our system \eqref{fn2} is only the first quadrant, the points $%
E_{11}$ and $E_{12}$ present interest and will be studied only when $\xi
_{11}\geq 0$ and $\xi _{12}\geq 0.$

\begin{definition}
We say an equilibrium point $E\left( \xi _{1},\xi _{2}\right) $ is proper if 
$\xi _{1}\geq 0$ and $\xi _{2}\geq 0,$ respectively, virtual if $\Delta \geq
0$ but $\xi _{1}<0$ or $\xi _{2}<0.$
\end{definition}

The first result is related to the type of bifurcation by which $E_{11}$ and 
$E_{12}$ come into existence or vanish.

\begin{proposition}
\label{prop-sn} If $\theta _{2}\delta \neq 0,$ $N\neq 0$ and $2N-\delta
\theta _{2}\neq 0,$ then $\Delta _{+}$ and $\Delta _{-}$ are saddle-node
bifurcation curves.
\end{proposition}

\textit{Proof.} Consider first the branch $\Delta _{+}.$ Write the system %
\eqref{fn2} in the form $\frac{d\xi }{dt}=f\left( \xi ,\mu \right) ,$ with $%
\xi =\left( \xi _{1},\xi _{2}\right) ,$ $f=\left( f_{1},f_{2}\right) $ and $%
\mu =\left( \mu _{1},\mu _{2}\right) .$ The proof will follow from
Sotomayor's theorem, as it is described in \cite{perko}. It is clear that $%
f\left( \xi _{0},\mu _{0}\right) =\left( 0,0\right) ,$ where $\xi
_{0}=\left( \xi _{11},0\right) $ and $\mu _{0}=\left( \mu _{1},\mu
_{2}\right) \in \Delta _{+}.$ Assume that $\mu _{2}>0$ is fixed while $\mu
_{1}$ varies, thus, $\mu _{1}$ is considered the bifurcation parameter.

The Jacobian matrices $A=Df\left( \xi _{0},\mu _{0}\right) =\left( 
\begin{array}{cc}
0 & \frac{\gamma \theta _{2}}{2N}\mu _{2}\left( 1+O\left( \mu _{2}\right)
\right)  \\ 
0 & \frac{2N-\delta \theta _{2}}{2N}\mu _{2}\left( 1+O\left( \mu _{2}\right)
\right) 
\end{array}%
\right) $ and $A^{T}$ have both an eigenvalue $\lambda =0,$ with the
corresponding eigenvectors $v$ for $A$ and $w$ for $A^{T},$ where $v=\left( 
\begin{array}{c}
1 \\ 
0%
\end{array}%
\right) \ $and $w=\left( 
\begin{array}{c}
\frac{\delta \theta _{2}-2N}{\gamma \theta _{2}}\left( 1+O\left( \mu
_{2}\right) \right)  \\ 
1%
\end{array}%
\right) .$ Denote by $f_{\mu _{1}}=\left( 
\begin{array}{c}
\frac{\partial f_{1}}{\partial \mu _{1}}\allowbreak  \\ 
\frac{\partial f_{2}}{\partial \mu _{1}}\allowbreak 
\end{array}%
\right) $ and $D^{2}f\left( \xi ,\mu \right) \left( v,v\right) =\left( 
\begin{array}{c}
d^{2}f_{1}\left( \xi ,\mu \right) \left( v,v\right)  \\ 
d^{2}f_{2}\left( \xi ,\mu \right) \left( v,v\right) 
\end{array}%
\right) ,$ where $d^{2}f_{1,2}\left( \xi ,\mu \right) \left( v,v\right) $
are the differentials of second order of the functions $f_{1,2}.$ Then 
\begin{equation*}
w^{T}f_{\mu _{1}}\left( \xi _{0},\mu _{0}\right) =-\mu _{2}\frac{2N-\delta
\theta _{2}}{2N\gamma }\left( 1+O\left( \mu _{2}\right) \right) \neq 0
\end{equation*}%
and%
\begin{equation*}
w^{T}\left[ D^{2}f\left( \xi _{0},\mu _{0}\right) \left( v,v\right) \right]
=-\mu _{2}\frac{2N-\delta \theta _{2}}{\gamma }\left( 1+O\left( \mu
_{2}\right) \right) \neq 0,
\end{equation*}%
which confirm the proof; if $w=\left( 
\begin{array}{c}
a \\ 
b%
\end{array}%
\right) $ then $w^{T}=\allowbreak \left( 
\begin{array}{cc}
a & b%
\end{array}%
\right) $ denotes the transpose vector of $w.$ For $\Delta _{-}$ the proof
is similar. $\blacksquare $

\begin{remark}
The first notable differences with the non-degenerate case \cite{GT1} is the
existence of two different equilibria $E_{11}$ and $E_{12}$ lying on the
same $\xi _{1}-$axis, and the existence of the saddle-node bifurcation
curves $\Delta _{+}$ and $\Delta _{-}.$ In the non-degenerate framework, a
single equilibrium $E_{1}$ exists on the $\xi _{1}-$axis and no saddle-node
bifurcation curves.
\end{remark}

\begin{remark}
$O(0,0)$ and $E_{2}\left( 0,-\mu _{2}+O\left( \mu _{2}^{2}\right) \right) $
are also equilibrium points of the system \eqref{fn2}. Their eigenvalues are 
$\mu _{1}$ and $\mu _{2}$ for $O,$ respectively, $\mu _{1}-\gamma \mu _{2}$
and $-\mu _{2}$ for $E_{2}\left( 0,-\mu _{2}\right) .$
\end{remark}

The system \eqref{fn2} has one more equilibrium point $E_{3}\left( \xi
_{1},\xi _{2}\right) ,$ where 
\begin{equation*}
\xi _{1}=-\frac{1}{\gamma \delta }\mu _{1}\left( 1+O\left( \left\vert \mu
\right\vert \right) \right) +\frac{1}{\delta }\mu _{2}\left( 1+O\left(
\left\vert \mu \right\vert \right) \right) \text{ and }\xi _{2}=-\frac{1}{%
\gamma }\mu _{1}\left( 1+O\left( \left\vert \mu \right\vert \right) \right)
+\sigma _{1}\mu _{2}^{2}\left( 1+O\left( \left\vert \mu \right\vert \right)
\right) ,
\end{equation*}%
with $\sigma _{1}=\frac{1}{\gamma \delta ^{2}}\left( \delta \theta
_{2}-N\right) .$ In its lowest terms, $E_{3}$ reads 
\begin{equation*}
E_{3}\left( -\frac{1}{\gamma \delta }\left( \mu _{1}-\gamma \mu _{2}\right)
,-\frac{1}{\gamma }\left( \mu _{1}-\gamma \sigma _{1}\mu _{2}^{2}\right)
\right) .
\end{equation*}

The existence and uniqueness of $E_{3}$ for $\left\vert \mu \right\vert $
sufficiently small is ensured by the Implicit Function Theorem applied to
the system$\qquad $%
\begin{equation}
\mu _{1}-\theta \left( \mu \right) \xi _{1}+\gamma \left( \mu \right) \xi
_{2}-M\left( \mu \right) \xi _{1}\xi _{2}+N\left( \mu \right) \xi _{1}^{2}=0%
\text{ and }\mu _{2}-\delta \left( \mu \right) \xi _{1}+\xi _{2}+S\left( \mu
\right) \xi _{1}^{2}+P\left( \mu \right) \xi _{2}^{2}=0.  \label{rel0}
\end{equation}

The bifurcation curves of $E_{3}$ for $\left\vert \mu \right\vert $
sufficiently small are

\begin{equation}
T_{2}=\left\{ \left( \mu _{1},\mu _{2}\right) \in 
%TCIMACRO{\U{211d} }%
%BeginExpansion
\mathbb{R}
%EndExpansion
^{2}\mid \mu _{2}=\frac{1}{\gamma }\mu _{1}+O\left( \mu _{1}^{2}\right) ,\mu
_{1}>0\right\}  \label{T2}
\end{equation}%
and

\begin{equation}
T_{3}=\left\{ \left( \mu _{1},\mu _{2}\right) \in 
%TCIMACRO{\U{211d} }%
%BeginExpansion
\mathbb{R}
%EndExpansion
^{2}\mid \mu _{1}=\gamma \sigma _{1}\mu _{2}^{2}+O\left( \mu _{2}^{3}\right)
,\delta \mu _{2}>0\right\} .  \label{T3}
\end{equation}

On $T_{2},$ $E_{3}$ coincides to $E_{2}\left( 0,-\mu _{2}\right) $ which
must have $\mu _{2}<0,$ while on $T_{3}$ to $E_{11}$ or $E_{12}.$ We call $%
E_{3}$ \textit{trivial} in these cases, otherwise nontrivial. For $%
\left\vert \mu \right\vert $ sufficiently small, $E_{3}$ is nontrivial in
the region 
\begin{equation}
R=\left\{ \left( \mu _{1},\mu _{2}\right) \in 
%TCIMACRO{\U{211d} }%
%BeginExpansion
\mathbb{R}
%EndExpansion
^{2}\mid \left( \mu _{1}-\gamma \mu _{2}\right) \delta >0,\mu _{1}-\gamma
\sigma _{1}\mu _{2}^{2}>0\right\} .  \label{r3}
\end{equation}

The characteristic polynomial at $E_{3}\left( \xi _{1},\xi _{2}\right) $ is $%
P\left( \lambda \right) =\lambda ^{2}-2p\lambda +L,$ where

\begin{equation}
p=\frac{1}{2}\left( \xi _{2}-\theta \left( \mu \right) \xi _{1}+2N\left( \mu
\right) \xi _{1}^{2}-M\left( \mu \right) \xi _{1}\xi _{2}+2P\left( \mu
\right) \xi _{2}^{2}\right)  \label{pin}
\end{equation}%
$\allowbreak $ and

\begin{equation}
L=\xi _{1}\xi _{2}\left( \gamma \left( \mu \right) \delta \left( \mu \right)
-\theta \left( \mu \right) +O\left( \left\vert \xi \right\vert \right)
\right) .  \label{q1}
\end{equation}%
For obtaining \eqref{pin} and \eqref{q1} we used \eqref{rel0}. The next
result describes the behavior of $E_{3}.$

\begin{theorem}
Assume $\theta _{2}\delta \neq 0,$ $N\neq 0$ and $\left( \mu _{1},\mu
_{2}\right) \in R.$ If $\delta >0,$ then $E_{3}$ is a saddle. If $\delta <0$
and $\sigma _{1}\neq 0,$ then

\qquad a) if $2N-\delta \theta _{2}>0,$ then $p>0$ and $E_{3}$ is a repeller;

\qquad b) if $2N-\delta \theta _{2}<0,$ then $E_{3}$ is an attractor (node
or focus) if $p<0$ and a repeller if $p>0.$ A Hopf bifurcation occurs at $%
E_{3}$ along the curve $p=0.$
\end{theorem}

Proof. Denoting by $\lambda _{1,2}=p\pm \sqrt{q}$ the eigenvalues of $E_{3},$
\eqref{q1} yields 
\begin{equation}
\lambda _{1}\lambda _{2}=\xi _{1}\xi _{2}\left( \gamma \delta +O\left(
\left\vert \mu \right\vert \right) \right) <0  \label{cp0}
\end{equation}%
if $\delta >0,$ thus, $E_{3}$ is a saddle; $\xi _{1,2}>0$ on $R$ and $\gamma
<0.$

Assume $\delta <0.$ This yields $\lambda _{1}\lambda _{2}>0$ by \eqref{q1}.
Then, $p$ in its lowest terms become%
\begin{equation*}
p=-\frac{1}{2\gamma }\mu _{1}+k_{3}\mu _{2}^{2}
\end{equation*}%
where $k_{3}=-\frac{1}{2\gamma \delta ^{2}}\left( N-\delta \theta
_{2}-\gamma \left( 2N-\delta \theta _{2}\right) \right) .$

1) Assume further $\sigma _{1}=\frac{1}{\gamma \delta ^{2}}\left( \delta
\theta _{2}-N\right) <0,$ which yields $N-\delta \theta _{2}<0$ and $%
T_{3}\subset \left\{ \mu _{1}>0,\mu _{2}<0\right\} .$ $R$ is included in the
fourth quadrant $\left\{ \mu _{1}>0,\mu _{2}<0\right\} $ from $\mu
_{1}>\gamma \sigma _{1}\mu _{2}^{2}>0$ and $\mu _{1}-\gamma \mu _{2}<0.$
Then $p=0$ occurs along the curve 
\begin{equation}
H=\left\{ \left( \mu _{1},\mu _{2}\right) \in 
%TCIMACRO{\U{211d} }%
%BeginExpansion
\mathbb{R}
%EndExpansion
^{2}\mid \mu _{1}=2\gamma k_{3}\mu _{2}^{2}+O\left( \mu _{2}^{3}\right) ,\mu
_{2}<0\right\} .  \label{h3}
\end{equation}

a) If $2N-\delta \theta _{2}>0$ then $\gamma \sigma _{1}>2\gamma k_{3}$ and $%
N>\delta \theta _{2}-N>0;$ $\gamma \sigma _{1}-2\gamma k_{3}=-\gamma \frac{%
2N-\delta \theta _{2}}{\delta ^{2}}.$ In this case, $p>0$ whenever $E_{3}$
exists, which yields that $E_{3}$ is a repeller.

b) If $2N-\delta \theta _{2}<0$ then $\gamma \sigma _{1}<2\gamma k_{3}$ and $%
H\subset R.$ Then, $p<0$ on the left of $H$ and $E_{3}$ is an attractor,
respectively, $p>0$ on the right of $H$ when $E_{3}$ is a repeller.

The eigenvalues of $E_{3}$ on $H$ are of the form $\pm i\omega _{0}$ where $%
\omega _{0}=\sqrt{\left. L\right\vert _{H}}=\frac{\mu _{2}}{\delta }\sqrt{%
-\gamma \mu _{2}\left( 2N-\delta \theta _{2}\right) }>0.$ Since $\left. 
\frac{\partial p}{\partial \mu _{1}}\right\vert _{H}=-\frac{1}{2\gamma }\neq
0,$ a Hopf bifurcation occurs on $H.$ It is nondegenerate if the first
Lyapunov coefficient $l_{1}\left( 0\right) \neq 0,$ otherwise, it is
degenerate.

2) Assume $\sigma _{1}>0,$ which yields $N-\delta \theta _{2}>0$ and $%
T_{3}\subset \left\{ \mu _{1}<0,\mu _{2}<0\right\} ;$ $\delta <0.$ In
addition, $R\subset \left\{ \mu _{2}<0\right\} .$

a) If $2N-\delta \theta _{2}>0$ which yields $\gamma \sigma _{1}-2\gamma
k_{3}=-\gamma \frac{2N-\delta \theta _{2}}{\delta ^{2}}>0,$ we get $%
H\nsubseteq R;$ $k_{3}>0.$ It follows that $p>0$ whenever $E_{3}$ exists,
that is, $E_{3}$ is a repeller.

b) If $2N-\delta \theta _{2}<0$ then $\gamma \sigma _{1}<2\gamma k_{3},$
which yields $H\subset R;$ in this case $N<\delta \theta _{2}-N<0.$ Similar
to 1b), $p<0$ on the left of $H$ and $E_{3}$ is an attractor, respectively, $%
p>0$ on the right of $H$ when $E_{3}$ is a repeller, either for $k_{3}>0$ or 
$k_{3}<0.$ The eigenvalues of $E_{3}$ on $H$ are of the form $\pm i\omega
_{0},$ thus, a Hopf bifurcation occurs on $H.$ $\blacksquare $

\begin{corollary}
Assume $\theta _{2}\delta \neq 0,$ $N>0$ and $\delta <0.$ If $\sigma _{1}<0$
then $\theta _{2}<0.$ If $\sigma _{1}>0$ then $2N-\delta \theta _{2}>0$ and $%
E_{3}$ is a repeller.
\end{corollary}

\begin{proposition}
\label{prop1} Assume $\theta _{2}\delta \neq 0,$ $N\neq 0$ and $\left( \mu
_{1},\mu _{2}\right) \in T_{3}.$ Then, in their lowest terms,

a) if $2N-\delta \theta _{2}>0,$ $E_{3}$ coincides to $E_{11}\left( \frac{1}{%
\delta }\mu _{2},0\right) ,$ $\lambda _{2}^{E_{11}}=0$ and $\lambda
_{1}^{E_{11}}>0.$ Moreover, $E_{12}\left( -\mu _{2}\frac{N-\delta \theta _{2}%
}{N\delta },0\right) $ is a saddle if $\mu _{2}N>0,$ respectively, an
attractor if $\mu _{2}N<0,$ whenever $E_{12}$ exists.

b) if $2N-\delta \theta _{2}<0,$ $E_{3}$ coincides to $E_{12}\left( \frac{1}{%
\delta }\mu _{2},0\right) ,$ $\lambda _{2}^{E_{12}}=0$ and $\lambda
_{1}^{E_{12}}<0.$ In addition, $E_{11}\left( -\mu _{2}\frac{N-\delta \theta
_{2}}{N\delta },0\right) $ is a saddle if $\mu _{2}N>0,$ respectively, a
repeller if $\mu _{2}N<0,$ whenever $E_{11}$ exists.
\end{proposition}

\textit{Proof}. Having the equilibrium point $E_{3}\left( \xi _{1},\xi
_{2}\right) ,$ $T_{3}$ is defined by $\xi _{2}=0.$ Other equilibria
satisfying $\xi _{2}=0$ are given by $\mu _{1}-\theta \left( \mu \right) \xi
_{1}+N\left( \mu \right) \xi _{1}^{2}=0,$ thus, they are $E_{11}$ or $%
E_{12}. $

The eigenvalues of an equilibrium point $\left( \xi _{1},0\right) $ are $%
3N\left( \mu \right) \xi _{1}^{2}-2\theta \left( \mu \right) \xi _{1}+\mu
_{1}$ and $S\left( \mu \right) \xi _{1}^{2}-\delta \left( \mu \right) \xi
_{1}+\mu _{2}.$ Thus, in their lowest terms, the eigenvalues of $%
E_{11}\left( \xi _{11},0\right) $ are 
\begin{equation}
\lambda _{1}^{E_{11}}=\xi _{11}\sqrt{\Delta }>0\text{ and }\lambda
_{2}^{E_{11}}=\frac{1}{N}\left( N\mu _{2}-S\mu _{1}-N\delta \xi
_{11}+S\theta _{2}\mu _{2}\xi _{11}\right) ,  \label{la10}
\end{equation}%
while of $E_{12}\left( \xi _{12},0\right) $ they read 
\begin{equation}
\lambda _{1}^{E_{12}}=-\xi _{12}\sqrt{\Delta }<0\text{ and }\lambda
_{2}^{E_{12}}=\frac{1}{N}\left( N\mu _{2}-S\mu _{1}-N\delta \xi
_{12}+S\theta _{2}\mu _{2}\xi _{12}\right) .  \label{la20}
\end{equation}

Let $\left( \mu _{1},\mu _{2}\right) \in T_{3},$ that is, $\mu _{1}=\gamma
\sigma _{1}\mu _{2}^{2}$ and $\delta \mu _{2}>0,$ which yield $\Delta
=\allowbreak \mu _{2}^{2}\frac{\left( 2N-\delta \theta _{2}\right) ^{2}}{%
\delta ^{2}}.$

a) If $2N-\delta \theta _{2}>0,$ then $E_{11}\left( \frac{1}{\delta }\mu
_{2},0\right) ,$ $E_{12}\left( -\mu _{2}\frac{N-\delta \theta _{2}}{N\delta }%
,0\right) $ and $E_{3}\left( \frac{1}{\delta }\mu _{2},0\right) ,$ thus, $%
E_{11}=E_{3}.$ Notice that $E_{3}\neq E_{12},$ otherwise, $2N-\delta \theta
_{2}=0.$ Therefore, $\lambda _{1}^{E_{3}}\lambda _{2}^{E_{3}}=\lambda
_{1}^{E_{11}}\lambda _{2}^{E_{11}}=0$ by \eqref{q1} since $\xi _{2}=0$ on $%
T_{3},$ where $\lambda _{1}^{E_{11}}=\xi _{11}\sqrt{\Delta }\neq 0,$ which
yield $\lambda _{2}^{E_{11}}=0$ and $\lambda _{1}^{E_{11}}>0,$ provided that 
$E_{11}$ is well-defined, that is, $\xi _{11}>0.$

Moreover, $\left. \lambda _{1}^{E_{12}}\right\vert _{T_{3}}=-\xi _{12}\sqrt{%
\Delta }<0$ and $\left. \lambda _{2}^{E_{12}}\right\vert _{T_{3}}=\mu _{2}%
\frac{2N-\delta \theta _{2}}{N}\left( 1+O\left( \mu _{2}\right) \right) >0$
if $\mu _{2}N>0,$ thus, $E_{12}$ is a saddle, respectively, $\left. \lambda
_{2}^{E_{12}}\right\vert _{T_{3}}<0$ if $\mu _{2}N<0,$ thus, $E_{12}$ is an
attractor, whenever $\xi _{12}>0.$

b) If $2N-\delta \theta _{2}<0,$ then $E_{11}\left( -\mu _{2}\frac{N-\delta
\theta _{2}}{N\delta },0\right) ,$ $E_{12}\left( \frac{1}{\delta }\mu
_{2},0\right) $ and $E_{3}\left( \frac{1}{\delta }\mu _{2},0\right) ,$ thus, 
$E_{3}=E_{12},$ $\lambda _{2}^{E_{12}}=0$ by \eqref{q1} and $\lambda
_{1}^{E_{12}}<0.$ In addition, $\lambda _{1}^{E_{11}}>0$ and $\lambda
_{2}^{E_{11}}\allowbreak =\mu _{2}\frac{2N-\delta \theta _{2}}{N}\left(
1+O\left( \mu _{2}\right) \right) ,$ whenever $\xi _{11}>0,$ thus, $E_{11}$
is a saddle if $\mu _{2}N>0,$ respectively, a repeller if $\mu _{2}N<0.$ $%
\blacksquare $

In the next theorem we characterize the bifurcation curve $T_{3}.$ Different
to the bifurcation curve $\Delta ,$ crossing $T_{3}$ in a non-degenerate
manner imposes a new condition, namely $\delta _{1}\neq 0.$

\begin{theorem}
\label{thT3} Assume $\theta _{2}\delta \delta _{1}\neq 0,$ $N\neq 0$ and $%
2N-\delta \theta _{2}\neq 0,$ where $\delta _{1}=\frac{\partial \delta }{%
\partial \mu _{1}}\left( 0\right) .$ Then $T_{3}$ is a transcritical
bifurcation curve.
\end{theorem}

\textit{Proof.} Assume $2N-\delta \theta _{2}>0.$ By Proposition \ref{prop1}%
, $E_{3}$ coincides to $E_{11}\left( \frac{1}{\delta }\mu _{2},0\right) $
and have the eigenvalues $\lambda _{2}^{E_{11}}=0$ and $\lambda
_{1}^{E_{11}}>0.$

Let $\mu _{1}$ be the bifurcation parameter while $\mu _{2}\neq 0$ is
assumed fixed. For $\mu _{2}\delta >0,$ denote by $\xi _{0}=\left( \frac{1}{%
\delta }\mu _{2},0\right) $ and $\mu _{0}=\left( \mu _{1},\mu _{2}\right)
\in T_{3}.$ Assume $\delta >0.$ The proof for $\delta <0$ is analogous.

Use further similar notations as in Proposition \ref{prop-sn}. Then $f\left(
\xi _{0},\mu _{0}\right) =\left( 0,0\right) ,$ respectively, $A=\left( 
\begin{array}{cc}
\frac{1}{\delta ^{2}}\mu _{2}^{2}\left( 2N-\delta \theta _{2}\right) & \frac{%
\gamma }{\delta }\mu _{2} \\ 
0 & 0%
\end{array}%
\right) ,$ $v=\left( 
\begin{array}{c}
v_{1} \\ 
1%
\end{array}%
\right) ,$ $w=\left( 
\begin{array}{c}
0 \\ 
1%
\end{array}%
\right) ,$ where $v_{1}=\frac{\gamma \delta }{\mu _{2}\left( \delta \theta
_{2}-2N\right) }.$ In finding $A,$ we imposed the conditions corresponding
to $E_{11},$ that is, $\xi _{2}=0,$ and to $E_{3},$ which read $\mu
_{1}-\theta \left( \mu \right) \xi _{1}+N\left( \mu \right) \xi _{1}^{2}=0$
and $\mu _{2}-\delta \left( \mu \right) \xi _{1}+S\left( \mu \right) \xi
_{1}^{2}=0;$ $A$ and $v_{1}$ are written in their lowest terms. Notice that $%
v_{1}$ is well-defined and $v_{1}\neq 0$ because $\left( 2N-\delta \theta
_{2}\right) \mu _{2}\neq 0$ and $0<\mu _{2}\ll 1.$

Denote by 
\begin{equation}
C_{1}=w^{T}f_{\mu _{1}}\left( \xi _{0},\mu _{0}\right) ,\text{ }C_{2}=w^{T}%
\left[ Df_{\mu _{1}}\left( \xi _{0},\mu _{0}\right) \left( v\right) \right] 
\text{ and }C_{3}=w^{T}\left[ D^{2}f\left( \xi _{0},\mu _{0}\right) \left(
v,v\right) \right] ,  \label{c123}
\end{equation}%
where $Df_{\mu _{1}}\left( \xi _{0},\mu _{0}\right) $ is the Jacobian matrix
in variables $\xi _{1}$ and $\xi _{2}$ of $f_{\mu _{1}}=\left( 
\begin{array}{cc}
\frac{\partial f_{1}}{\partial \mu _{1}} & \frac{\partial f_{2}}{\partial
\mu _{1}}%
\end{array}%
\right) ^{T}$ calculated at $\left( \xi _{0},\mu _{0}\right) .$

In order to determine $f_{\mu _{1}}$ from the system \eqref{fn2}, we need to
write 
\begin{equation*}
\delta \left( \mu \right) =\delta \left( 0\right) +\delta _{1}\mu _{1}\left(
1+O\left( \mu _{1}\right) \right) +\delta _{2}\mu _{2}\left( 1+O\left( \mu
_{2}\right) \right) ,
\end{equation*}%
and similarly the other parameter-functions $\theta \left( \mu \right) ,$ $%
\gamma \left( \mu \right) ,$ $M\left( \mu \right) $ and so on; $\delta _{1}=%
\frac{\partial \delta }{\partial \mu _{1}}\left( 0\right) $ and $\delta _{2}=%
\frac{\partial \delta }{\partial \mu _{2}}\left( 0\right) .$ However, only $%
\delta \left( \mu \right) $ will be needed in this case.

Then, $f_{\mu _{1}}$ in its lowest terms in $\mu _{1}$ and $\mu _{2}$ has
the form 
\begin{equation}
f_{\mu _{1}}=\left( 
\begin{array}{c}
\xi _{1}\left( 1-\theta _{1}\xi _{1}+\gamma _{1}\xi _{2}-M_{1}\xi _{1}\xi
_{2}+N_{1}\xi _{1}^{2}\right) \\ 
\xi _{2}\left( -\delta _{1}\xi _{1}+S_{1}\xi _{1}^{2}+P_{1}\xi
_{2}^{2}\right)%
\end{array}%
\right) ,  \label{fmu1}
\end{equation}%
where $\theta _{1}=\frac{\partial \theta }{\partial \mu _{1}}\left( 0\right)
,$ $\gamma _{1}=\frac{\partial \gamma }{\partial \mu _{1}}\left( 0\right) ,$ 
$M_{1}=\frac{\partial M}{\partial \mu _{1}}\left( 0\right) $ and so on. By %
\eqref{c123}, these lead to

\begin{equation*}
C_{1}=0,\text{ }C_{2}=-\frac{\delta _{1}}{\delta }\mu _{2}\left( 1+O\left(
\mu _{2}\right) \right) \neq 0\text{ and }C_{3}=\frac{2\gamma \delta ^{2}}{%
\left( 2N-\delta \theta _{2}\right) \mu _{2}}\left( 1+O\left( \mu
_{2}\right) \right) \neq 0,
\end{equation*}%
because $0<\mu _{2}\ll 1$ and $2N-\delta \theta _{2}>0.$ The proof is
similar for $2N-\delta \theta _{2}<0.$ Thus, a non-degenerate transcritical
bifurcation occurs on $T_{3}.\blacksquare $

\bigskip

Denote by 
\begin{equation*}
X_{+}=\left\{ \left( \mu _{1},\mu _{2}\right) \left\vert \mu _{1}>0,\mu
_{2}=0\right. \right\} ,\text{ }X_{-}=\left\{ \left( \mu _{1},\mu
_{2}\right) \left\vert \mu _{1}<0,\mu _{2}=0\right. \right\} ,
\end{equation*}%
respectively, 
\begin{equation*}
Y_{+}=\left\{ \left( \mu _{1},\mu _{2}\right) \left\vert \mu _{2}>0,\mu
_{1}=0\right. \right\} ,\text{ }Y_{-}=\left\{ \left( \mu _{1},\mu
_{2}\right) \left\vert \mu _{2}<0,\mu _{1}=0\right. \right\} ,
\end{equation*}%
the semi-major and semi-minor axes of coordinates. The next result describes
the bifurcations which occur on the remaining curves. Different to $T_{3},$
they do not need the constraint $\delta _{1}\neq 0.$

\begin{theorem}
\label{thT4} Assume $\theta _{2}\delta \neq 0,$ $N\neq 0$ and $2N-\delta
\theta _{2}\neq 0.$ Then $T_{2},$ $Y_{+},$ $Y_{-},$ $X_{+}$ and $X_{-}$ are
transcritical bifurcation curves.
\end{theorem}

\textit{Proof.} When $\left( \mu _{1},\mu _{2}\right) \in T_{2},$ $%
E_{2}\left( 0,-\mu _{2}+O\left( \mu _{2}^{2}\right) \right) $ coincides to $%
E_{3},$ thus, $\xi _{0}=\left( 0,-\mu _{2}+O\left( \mu _{2}^{2}\right)
\right) $ and $\mu _{0}=\left( \mu _{1},\mu _{2}\right) \in T_{2}$ with $\mu
_{2}\neq 0.$ We find $v=\left( 
\begin{array}{cc}
\frac{1}{\delta }\left( 1+O\left( \mu _{2}\right) \right) & 1%
\end{array}%
\right) ^{T}$ and $w=\left( 
\begin{array}{cc}
1 & 0%
\end{array}%
\right) ^{T}.$ Using $\mu _{1}$ as the bifurcation parameter, \eqref{fmu1}
and \eqref{c123} yield

\begin{equation*}
C_{1}=0,\text{ }C_{2}=\frac{1}{\delta }\left( 1+O\left( \mu _{2}\right)
\right) \neq 0\text{ and }C_{3}=\frac{2\gamma }{\delta }\left( 1+O\left( \mu
_{2}\right) \right) \neq 0,
\end{equation*}%
thus, the bifurcation on $T_{2}$ is transcritical.

Let $\left( 0,\mu _{2}\right) \in Y_{+}\cup Y_{-}.$ Then $O$ coincides to $%
E_{12}$ if $\theta _{2}\mu _{2}>0,$ respectively, $E_{11}$ if $\theta
_{2}\mu _{2}<0.$ Also, $\xi _{0}=\left( 0,0\right) ,$ $\mu _{0}=\left( 0,\mu
_{2}\right) $ with $\mu _{2}\neq 0,$ and $v=w=\allowbreak \left( 
\begin{array}{cc}
1 & 0%
\end{array}%
\right) ^{T};$ the bifurcation parameter is $\mu _{1}.$ These yield $%
C_{1}=0, $ $C_{2}=1$ and $C_{3}=-2\theta _{2}\mu _{2}\neq 0.$

Let $\left( \mu _{1},0\right) \in X_{+}\cup X_{-}$ and consider $\mu _{2}$
as the bifurcation parameter. Then $O$ coincides to $E_{2}$ and $\xi
_{0}=\left( 0,0\right) ,$ $\mu _{0}=\left( \mu _{1},0\right) $ with $\mu
_{1}\neq 0,$ and $v=w=\allowbreak \left( 
\begin{array}{cc}
0 & 1%
\end{array}%
\right) ^{T}.$ These lead to $w^{T}f_{\mu _{2}}\left( \xi _{0},\mu
_{0}\right) =0,$ $w^{T}\left[ Df_{\mu _{2}}\left( \xi _{0},\mu _{0}\right)
\left( v\right) \right] =1$ and $w^{T}\left[ D^{2}f\left( \xi _{0},\mu
_{0}\right) \left( v,v\right) \right] =2.$ $\blacksquare $

The next result is important because it states that the signs of the
eigenvalues of $E_{11}$ and $E_{12}$ depend on the conditions of existence
of $E_{3}.$

\begin{proposition}
\label{prop3} Assume $\theta _{2}\delta \neq 0,$ $N\neq 0$ and $\delta \mu
_{2}>0.$ If $2N-\delta \theta _{2}>0,$ the curve $\lambda _{2}^{E_{11}}=0$
is unique and coincides to $T_{3}$ for $\left\vert \mu \right\vert $
sufficiently small. Similarly, if $2N-\delta \theta _{2}<0,$ $\lambda
_{2}^{E_{12}}=0$ coincides to $T_{3}.$
\end{proposition}

\textit{Proof.} By \eqref{la10} - \eqref{la20} we have

\begin{equation}
\lambda _{2}^{E_{11}}\lambda _{2}^{E_{12}}=\frac{\delta ^{2}}{N}\left[ \mu
_{1}\left( 1+O\left( \left\vert \mu \right\vert \right) \right) -\gamma
\sigma _{1}\mu _{2}^{2}\left( 1+O\left( \left\vert \mu \right\vert \right)
\right) \right] ,  \label{la22}
\end{equation}%
where $\sigma _{1}=\frac{1}{\gamma \delta ^{2}}\left( \delta \theta
_{2}-N\right) .$ From the Implicit Function Theorem applied to the
right-hand side term of \eqref{la22}, there exists a unique curve $%
T_{3}^{\prime }$ such that $\lambda _{2}^{E_{11}}\lambda _{2}^{E_{12}}=0$ on 
$T_{3}^{\prime },$ for $\left\vert \mu \right\vert $ sufficiently small and $%
\delta \mu _{2}>0,$ given by

\begin{equation}
T_{3}^{\prime }=\left\{ \left( \mu _{1},\mu _{2}\right) \in 
%TCIMACRO{\U{211d} }%
%BeginExpansion
\mathbb{R}
%EndExpansion
^{2}\mid \mu _{1}=\gamma \sigma _{1}\mu _{2}^{2}+k_{1}\mu _{2}^{3}+O\left(
\mu _{2}^{4}\right) ,\delta \mu _{2}>0\right\} ,  \label{t3prim}
\end{equation}%
where $k_{1}\in 
%TCIMACRO{\U{211d} }%
%BeginExpansion
\mathbb{R}
%EndExpansion
.$ Therefore, $\lambda _{2}^{E_{11}}=0$ or $\lambda _{2}^{E_{12}}=0$ on $%
T_{3}^{\prime }.$ But $\left. \lambda _{2}^{E_{12}}\right\vert
_{T_{3}^{\prime }}=\frac{2N-\delta \theta _{2}}{N}\mu _{2}\left( 1+O\left(
\mu _{2}\right) \right) \neq 0$ for $\left\vert \mu \right\vert $
sufficiently small, thus, $\left. \lambda _{2}^{E_{11}}\right\vert
_{T_{3}^{\prime }}=0.$ On the other hand, from Proposition \ref{prop1} we
know also $\left. \lambda _{2}^{E_{11}}\right\vert _{T_{3}}=0.$ But $%
T_{3}^{\prime }$ is unique from the Implicit Function Theorem, thus $%
T_{3}=T_{3}^{\prime }$ for $\left\vert \mu \right\vert $ sufficiently small. 
$\blacksquare $

\bigskip

\textbf{Assume further }$N>0$ and $\theta _{1}=\frac{\partial \theta }{%
\partial \mu _{1}}\left( 0\right) \neq 0.$ The case $N<0$ can be treated
similarly. Define the following regions 
\begin{equation*}
R_{00}=\left\{ \left( \mu _{1},\mu _{2}\right) \left\vert \Delta <0\right.
\right\} \cup \left\{ \left( \mu _{1},\mu _{2}\right) \left\vert \Delta
>0,\mu _{1}>0,\theta _{2}\mu _{2}<0\right. \right\} ,
\end{equation*}

\begin{equation}
R_{10}=\left\{ \left( \mu _{1},\mu _{2}\right) \left\vert \mu _{1}<0\right.
\right\} \text{ and }R_{20}=\left\{ \left( \mu _{1},\mu _{2}\right)
\left\vert \Delta >0,\mu _{1}>0,\theta _{2}\mu _{2}>0\right. \right\} .
\label{r10}
\end{equation}

Since $\xi _{11}\xi _{12}=\mu _{1}/N,$ it follows that on $R_{10}$ a single
proper equilibrium exists, $E_{11},$ while $E_{12}$ is virtual. On $R_{20},$
both equilibria $E_{11}$ and $E_{12}$ exist. Indeed, in their lowest terms,
we have 
\begin{equation*}
\xi _{11}+\xi _{12}=\frac{1}{N}\theta \left( \mu \right) =\frac{1}{N}\left(
\theta _{1}\mu _{1}+\theta _{2}\mu _{2}\right) .
\end{equation*}%
If $\theta _{1}>0,$ it is clear that $\xi _{11}>0$ and $\xi _{12}>0$ on $%
R_{20}.$ Assume $\theta _{1}<0$ and $\left( \mu _{1},\mu _{2}\right) \in
R_{20}.$ Then $\mu _{1}\theta _{1}>\frac{\theta _{2}^{2}}{4N}\mu
_{2}^{2}\theta _{1}$ from $\Delta >0,$ which leads to $\theta _{1}\mu
_{1}+\theta _{2}\mu _{2}>\allowbreak \theta _{2}\mu _{2}\left( 1+\frac{%
\theta _{2}\theta _{1}}{4N}\mu _{2}\right) >0$ for $\left\vert \mu
\right\vert $ sufficiently small. Thus, $\xi _{11}>0$ and $\xi _{12}>0$ on $%
R_{20}.$

One can show similarly that on $R_{00},$ $E_{11}$ and $E_{12}$ do not exist
because either $\xi _{11}$ and $\xi _{12}$ are not real numbers ($\Delta <0$%
) or $\xi _{11}<0$ and $\xi _{12}<0$ ($E_{11}$ and $E_{12}$ are virtual
points) because $\theta _{1}\mu _{1}+\theta _{2}\mu _{2}<\allowbreak \theta
_{2}\mu _{2}\left( 1+\frac{\theta _{2}\theta _{1}}{4N}\mu _{2}\right) <0.$

Whenever $\sigma _{1}<0,$ denote by $R_{20}^{-}=R_{20}\cap \left\{ \left(
\mu _{1},\mu _{2}\right) \left\vert \mu _{1}<\gamma \sigma _{1}\mu
_{2}^{2}\right. \right\} $ and $R_{20}^{+}=R_{20}\cap \left\{ \left( \mu
_{1},\mu _{2}\right) \left\vert \mu _{1}>\gamma \sigma _{1}\mu
_{2}^{2}\right. \right\} ,$ the regions from $R_{20}$ to the left,
respectively, the right of $T_{3}.$ Notice that 
\begin{equation*}
R_{20}=R_{20}^{-}\cup T_{3}\cup R_{20}^{+}.
\end{equation*}

\begin{theorem}
\label{th5} Assume $N>0$ and $\sigma _{1}<0.$ If $\delta >0,$ then

a) if $2N-\delta \theta _{2}>0,$ then $E_{11}$ is a saddle on $R_{10}\cup
Y_{+}\cup R_{20}^{-}$ and a repeller on $R_{20}^{+},$ while $E_{12}$ is a
saddle on $R_{20};$

b) if $2N-\delta \theta _{2}<0,$ then $E_{11}$ is a saddle on $R_{10}\cup
Y_{+}\cup R_{20},$ while $E_{12}$ is a saddle on $R_{20}^{-}$ and an
attractor on $R_{20}^{+}.$
\end{theorem}

\textit{Proof.} With the help of \eqref{la22}, we are able to determine the
dynamics of $E_{11}$ and $E_{12}$ when $\left\vert \mu \right\vert $ is
sufficiently small.

From $N>0,$ $\sigma _{1}<0$ and $\delta >0$ we get $N-\delta \theta _{2}<0,$ 
$\theta _{2}>0$ and $\frac{\theta _{2}^{2}}{4N}-\gamma \sigma _{1}=\frac{%
\left( 2N-\delta \theta _{2}\right) ^{2}}{4N\delta ^{2}}>0;$ assume $%
2N-\delta \theta _{2}\neq 0.$ $R_{20}$ and $\Delta _{+}$ lie on $\left\{ \mu
_{1}>0,\mu _{2}>0\right\} .$

We describe in the following the behavior of the points $E_{11}$ and $%
E_{12}. $ When $\left( \mu _{1},\mu _{2}\right) $ crosses $\Delta _{+},$ the
points $E_{11}$ and $E_{12}$ are born in the region $\Delta >0$ through a
saddle-node bifurcation on the curve $\Delta _{+}.$

On $\Delta _{+},$ there is a single equilibrium of the form $\left( \xi
_{1},0\right) ,$ namely $E_{11}\left( \frac{\theta _{2}}{2N}\mu
_{2},0\right) $ which has the eigenvalues $0$ and $\lambda _{2}^{E_{11}}=\mu
_{2}\frac{2N-\delta \theta _{2}}{2N}\left( 1+O\left( \mu _{2}\right) \right)
;$ $E_{11}$ coincides to $E_{12}$ on $\Delta _{+}$ while $E_{3}\left( \frac{%
\mu _{2}}{\delta },-\mu _{2}^{2}\frac{\left( 2N-\delta \theta _{2}\right)
^{2}}{4\gamma \delta ^{2}N}\right) $ may exist as a different point.

As soon as the point $\left( \mu _{1},\mu _{2}\right) $ leaves the curve $%
\Delta _{+}$ and $\left( \mu _{1},\mu _{2}\right) \in R_{20},$ $E_{11}$ and $%
E_{12}$ exist as two different points. In order to study the behavior of $%
E_{11}$ and $E_{12}$ when $\left( \mu _{1},\mu _{2}\right) \in R_{20},$ we
will use \eqref{la22}. $\lambda _{2}^{E_{11}}\lambda _{2}^{E_{12}}=0$ on $%
T_{3}$ yields $\lambda _{2}^{E_{11}}=0$ or $\lambda _{2}^{E_{12}}=0$ on $%
T_{3};$ they cannot be at the same time $0$ on $T_{3}$ because $\lambda
_{2}^{E_{11}}-\lambda _{2}^{E_{12}}=-\frac{1}{N}\left( N\delta -S\theta
_{2}\mu _{2}\right) \left( \xi _{11}-\xi _{12}\right) \neq 0$ if $\Delta >0$
and $\left\vert \mu \right\vert $ sufficiently small.

a) Assume $2N-\delta \theta _{2}>0$ and $\left( \mu _{1},\mu _{2}\right) \in
T_{3},$ i.e. $\mu _{1}=\gamma \sigma _{1}\mu _{2}^{2}$ and $\delta \mu
_{2}>0.$ Then $\xi _{12}=-\mu _{2}\frac{N-\delta \theta _{2}}{N\delta }>0$
and 
\begin{equation}
\lambda _{2}^{E_{12}}=\mu _{2}\frac{2N-\delta \theta _{2}}{N}\left(
1+O\left( \mu _{2}\right) \right) \neq 0  \label{la2e12}
\end{equation}%
on $T_{3};$ $\mu _{2}>0$ on $T_{3}$ since $\delta \mu _{2}>0$ by \eqref{T3}.
Therefore, $\lambda _{2}^{E_{11}}=0$ on $T_{3}$ and, thus, $\lambda
_{2}^{E_{11}}$ changes its sign when $\left( \mu _{1},\mu _{2}\right) $
crosses $T_{3},$ while $\lambda _{2}^{E_{12}}\neq 0$ on $T_{3}$ and, thus, $%
\lambda _{2}^{E_{12}}$ keeps constant sign on $\mu _{2}>0$ sufficiently
small, namely $\lambda _{2}^{E_{12}}>0$ by \eqref{la2e12}.

By \eqref{la22}, $\lambda _{2}^{E_{11}}\lambda _{2}^{E_{12}}>0$ on the right
of $T_{3},$ i.e. on $R_{20}^{+}.$ Thus, $\lambda _{2}^{E_{11}}>0$ on $%
R_{20}^{+}$ and $\lambda _{2}^{E_{11}}<0$ on the left of $T_{3},$ i.e. on $%
R_{20}^{-}.$ Using $\lambda _{1}^{E_{11}}>0$ and $\lambda _{1}^{E_{12}}<0$
whenever $E_{11}$ and $E_{12}$ exist, it follows that $E_{11}$ is a repeller
on $R_{20}^{+}$ and a saddle on $R_{20}^{-},$ while $E_{12}$ is a saddle on $%
R_{20}.$

On $T_{3},$ $E_{11}\left( \frac{1}{\delta }\mu _{2},0\right) $ having the
eigenvalues $\left. \lambda _{2}^{E_{11}}\right\vert _{T_{3}}=0$ and $\left.
\lambda _{1}^{E_{11}}\right\vert _{T_{3}}>0,$ coincides to $E_{3}.$ $E_{12}$
is a saddle on $T_{3},$ because $\left. \lambda _{2}^{E_{12}}\right\vert
_{T_{3}}=\mu _{2}\frac{2N-\delta \theta _{2}}{N}\left( 1+O(\mu _{2})\right)
>0$ for $\mu _{2}>0$ sufficiently small and $\left. \lambda
_{1}^{E_{12}}\right\vert _{T_{3}}<0.$

On $Y_{+},$ $E_{11}\left( \frac{\theta _{2}}{N}\mu _{2},0\right) $ continues
to survive as a saddle point while $E_{12}$ collides to $O;$ $\lambda
_{2}^{E_{11}}=\frac{N-\delta \theta _{2}}{N}\mu _{2}\left( 1+O(\mu
_{2})\right) <0$ on $Y_{+}.$ On $Y_{-},$ $E_{11}$ collides to $O$ and
vanishes on $\mu _{1}>0$ and $\mu _{2}>0.$

On $\mu _{1}<0,$ $\lambda _{2}^{E_{11}}$ keeps constant (negative) sign, $%
\lambda _{2}^{E_{11}}<0,$ because $T_{3}\nsubseteq \left\{ \mu
_{1}<0\right\} .$ This is in agreement with 
\begin{equation}
\lambda _{2}^{E_{11}}\left( \mu _{1},0\right) =\frac{\sqrt{-N\mu _{1}}}{N}%
\left( \frac{S}{N}\sqrt{-N\mu _{1}}-\delta \right) <0  \label{la02}
\end{equation}%
for $\mu _{1}<0$ sufficiently small; $\delta >0$ and $\xi _{11}\left( \mu
_{1},0\right) =\frac{1}{2N}\sqrt{-4N\mu _{1}}>0$ on $\mu _{1}<0.$ It implies
that $E_{11}$ survives in $R_{10}$ as a saddle point while $E_{12}$ vanishes
in $R_{10}.$

On $\Delta _{+},$ the eigenvalues of the coinciding points $E_{11}\left( 
\frac{\theta _{2}}{2N}\mu _{2},0\right) $ and $E_{12}$ are in this case $0$
and $\lambda _{2}^{E_{11}}=\mu _{2}\frac{2N-\delta \theta _{2}}{2N}\left(
1+O\left( \mu _{2}\right) \right) >0.$

b) Assume $2N-\delta \theta _{2}<0.$ Then $\xi _{11}=-\mu _{2}\frac{N-\delta
\theta _{2}}{N\delta }>0$ and 
\begin{equation}
\lambda _{2}^{E_{11}}=\mu _{2}\frac{2N-\delta \theta _{2}}{N}\left(
1+O\left( \mu _{2}\right) \right) <0  \label{la2e11}
\end{equation}%
on $T_{3}.$ Therefore, $\lambda _{2}^{E_{12}}=0$ on $T_{3}$ and $\lambda
_{2}^{E_{12}}$ changes its sign when $\left( \mu _{1},\mu _{2}\right) $
crosses $T_{3},$ while $\lambda _{2}^{E_{11}}$ keeps constant (negative)
sign on $R_{10}\cup Y_{+}\cup R_{20}.$ Thus, $E_{11}$ is a saddle on $%
R_{10}\cup Y_{+}\cup R_{20},$ including on $T_{3}\subset R_{20}.$

From $\lambda _{2}^{E_{11}}\lambda _{2}^{E_{12}}>0$ on $R_{20}^{+},$ we have 
$\lambda _{2}^{E_{12}}<0$ on $R_{20}^{+}$ and $\lambda _{2}^{E_{12}}>0$ on $%
R_{20}^{-}.$ Therefore, $E_{12}$ is an attractor on $R_{20}^{+}$ and a
saddle on $R_{20}^{-}$ because $\lambda _{1}^{E_{12}}<0.$ On $Y_{+}\cup
Y_{-} $ the results are similar to a).

On $T_{3},$ $E_{3}$ coincides to $E_{12}\left( \frac{1}{\delta }\mu
_{2},0\right) ,$ which has $\left. \lambda _{2}^{E_{12}}\right\vert
_{T_{3}}=0$ and $\left. \lambda _{1}^{E_{12}}\right\vert _{T_{3}}<0.$ On $%
\Delta _{+},$ the eigenvalues of $E_{11}$ are $0$ and $\lambda
_{2}^{E_{11}}=\mu _{2}\frac{2N-\delta \theta _{2}}{2N}\left( 1+O\left( \mu
_{2}\right) \right) <0.$ $\blacksquare $

\begin{theorem}
Assume $N>0$ and $\sigma _{1}<0.$ If $\delta <0,$ then

a) if $2N-\delta \theta _{2}>0,$ $E_{11}$ is a repeller on $R_{10}\cup
Y_{-}\cup R_{20}^{-}$ and a saddle on $R_{20}^{+},$ while $E_{12}$ is an
attractor on $R_{20}.$

b) if $2N-\delta \theta _{2}<0,$ $E_{11}$ is a repeller on $R_{10}\cup
Y_{-}\cup R_{20},$ while $E_{12}$ is an attractor on $R_{20}^{-}$ and a
saddle on $R_{20}^{+}.$
\end{theorem}

\textit{Proof.} The hypothesis leads to $N-\delta \theta _{2}<0,$ $\theta
_{2}<0$ and $\frac{\theta _{2}^{2}}{4N}>\gamma \sigma _{1},$ where $%
2N-\delta \theta _{2}\neq 0.$ Different to the first case, now $T_{3}$ and $%
R_{20}$ lie in the fourth quadrant, $\left\{ \mu _{1}>0,\mu _{2}<0\right\} ,$
and the used branch of $\Delta $ is $\Delta _{-}.$

a) Let $2N-\delta \theta _{2}>0.$ Then $\lambda _{2}^{E_{12}}\neq 0$ on $%
T_{3}$ by \eqref{la2e12} and, thus, $\lambda _{2}^{E_{11}}=0$ on $T_{3}$ by %
\eqref{la22}. Therefore, $\lambda _{2}^{E_{12}}<0$ on $R_{20}$ while $%
\lambda _{2}^{E_{11}}$ changes its sign when crossing $T_{3}.$ More exactly,
from $\lambda _{2}^{E_{11}}\lambda _{2}^{E_{12}}>0$ on $R_{20}^{+}$ by %
\eqref{la22}, we have $\lambda _{2}^{E_{11}}<0$ on $R_{20}^{+}$ and $\lambda
_{2}^{E_{11}}>0$ on $R_{20}^{-}.$ Thus, $E_{11}$ is a repeller on $%
R_{20}^{-} $ and a saddle on $R_{20}^{+},$ respectively, $E_{12}$ is an
attractor on $R_{20}.$

On $Y_{+},$ $E_{11}$ collides to $O$ while $E_{12}$ exists only virtually
because $\xi _{12}<0.$ However, as soon as $\mu _{1}<0,$ $E_{11}$ bifurcates
from $O.$

On $R_{10},$ $E_{11}$ is a repeller; $\lambda _{2}^{E_{11}}>0$ on $R_{10}$
follows from \eqref{la02} and $T_{3}\nsubseteq \left\{ \mu _{1}<0\right\} .$

On $Y_{-},$ $E_{11}\left( \frac{\theta _{2}}{N}\mu _{2},0\right) $ remains a
repeller since $\lambda _{1}^{E_{11}}=\xi _{11}\sqrt{\Delta }>0$ and $%
\lambda _{2}^{E_{11}}=\frac{N-\delta \theta _{2}}{N}\mu _{2}\left( 1+O\left(
\mu _{2}\right) \right) >0,$ while $E_{12}$ collides to $O.$

On $\Delta _{-},$ the eigenvalues of the coinciding points $E_{11}$ and $%
E_{12}$ are $0$ and $\lambda _{2}^{E_{11}}=\mu _{2}\frac{2N-\delta \theta
_{2}}{2N}\left( 1+O\left( \mu _{2}\right) \right) <0.$

b) Let $2N-\delta \theta _{2}<0.$ Then $\lambda _{2}^{E_{11}}\neq 0$ and $%
\lambda _{2}^{E_{12}}=0$ on $T_{3}$ by \eqref{la2e11}. From $\xi
_{11}=\allowbreak -\mu _{2}\frac{N-\delta \theta _{2}}{N\delta }>0$ and $%
\lambda _{2}^{E_{11}}=\mu _{2}\frac{2N-\delta \theta _{2}}{N}\left(
1+O\left( \mu _{2}\right) \right) >0$ on $T_{3},$ it follows that $\lambda
_{2}^{E_{11}}>0$ whenever $E_{11}$ exists. Thus, $E_{11}$ is a repeller on $%
R_{10}\cup Y_{-}\cup R_{20}.$ On $Y_{+},$ $E_{11}$ collides to $O.$

Further, $\lambda _{2}^{E_{11}}\lambda _{2}^{E_{12}}>0$ on $R_{20}^{+}$
yields $\lambda _{2}^{E_{12}}>0$ on $R_{20}^{+}$ and $\lambda
_{2}^{E_{12}}<0 $ on $R_{20}^{-}.$ Thus, $E_{12}$ is a saddle on $R_{20}^{+}$
and an attractor on $R_{20}^{-}.$

On $\Delta _{-},$ the eigenvalue $\lambda _{2}^{E_{11}}$ becomes positive, $%
\lambda _{2}^{E_{11}}=\mu _{2}\frac{2N-\delta \theta _{2}}{2N}\left(
1+O\left( \mu _{2}\right) \right) >0,$ while the other remains $0.$ $%
\blacksquare $

\begin{figure}[tbp]
\centering
\includegraphics[width=1\textwidth]{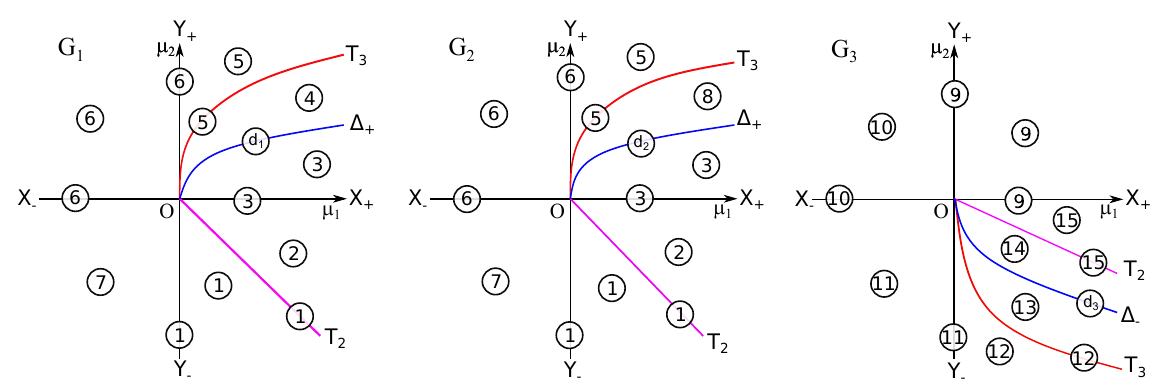} \includegraphics[width=1%
\textwidth]{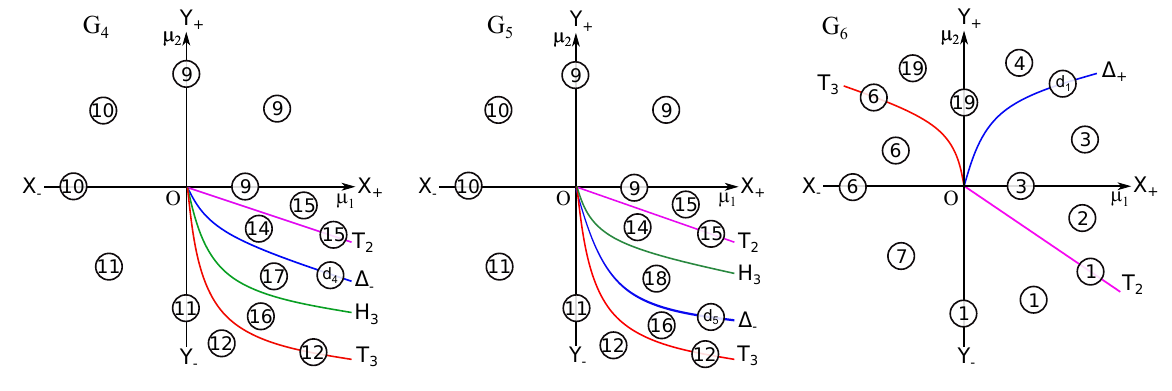}
\caption{Bifurcation diagrams corresponding to $N>0,\protect\sigma _{1}<0$
and (G1) $\protect\delta >0,$ $2N-\protect\delta \protect\theta _{1}>0,$
(G2) $\protect\delta >0,$ $2N-\protect\delta \protect\theta _{1}<0,$ (G3) $%
\protect\delta <0,$ $2N-\protect\delta \protect\theta _{1}>0,$ (G4) $\protect%
\delta <0,$ $2N-\protect\delta \protect\theta _{1}<0,$ $2\protect\gamma k_3<%
\frac{\protect\theta_1^2}{4N},$ and (G5) $\protect\delta <0,$ $2N-\protect%
\delta \protect\theta _{1}<0,$ $2\protect\gamma k_3>\frac{\protect\theta_1^2%
}{4N},$ respectively, (G6) corresponding to $N>0,$ $\protect\sigma_1>0,$ $%
\protect\delta>0$ and $\protect\theta_1>0.$}
\label{pp3}
\end{figure}

\begin{remark}
When $N>0,$ $\sigma _{1}<0,$ $\delta <0$ and $2N-\delta \theta _{2}<0,$ the
curve $H$ can lie on the both sides of the curve $\Delta _{-}.$ More
exactly, $H$ lies on the left of $\Delta _{-}$ if $2N-4N\gamma -\delta
\theta _{2}<0,$ respectively, the right of $\Delta _{-}$ if $2N-4N\gamma
-\delta \theta _{2}>0,$ because $2\gamma k_{3}-\frac{\theta _{2}^{2}}{4N}%
=\allowbreak \frac{\delta \theta _{2}-2N}{4N\delta ^{2}}\left( 2N-4N\gamma
-\delta \theta _{2}\right) ,$ Fig. \ref{pp3}.
\end{remark}

\begin{figure}[tbp]
\centering
\includegraphics[width=1\textwidth]{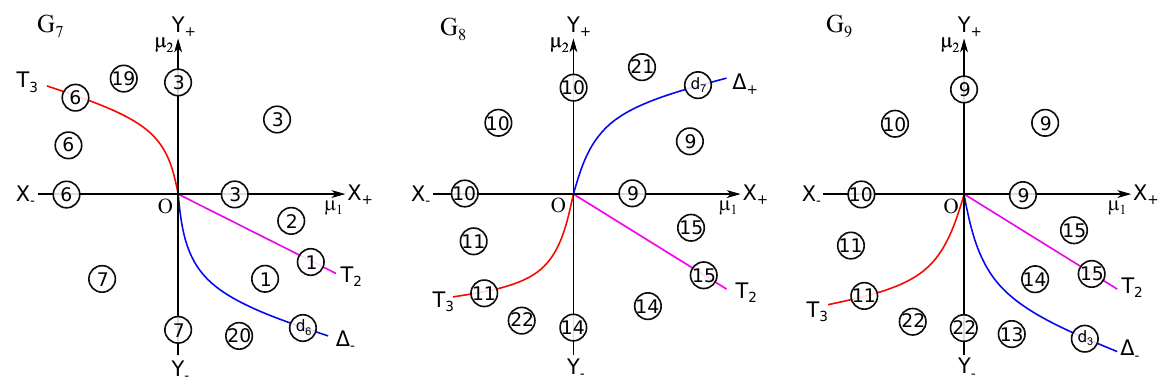}
\caption{Bifurcation diagrams corresponding to $N>0,\protect\sigma_1>0$ and
(G7) $\protect\delta>0,$ $\protect\theta_1<0,$ (G8) $\protect\delta<0,$ $%
\protect\theta_1>0,$ and (G9) $\protect\delta<0,$ $\protect\theta_1<0.$}
\label{pp4}
\end{figure}

\begin{figure}[tbp]
\centering
\includegraphics[width=0.45\textwidth]{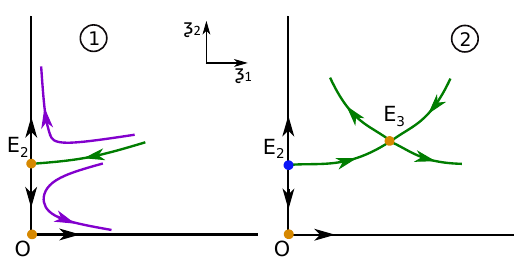} %
\includegraphics[width=0.45\textwidth]{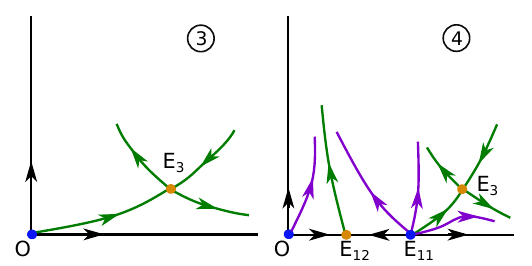} %
\includegraphics[width=0.45\textwidth]{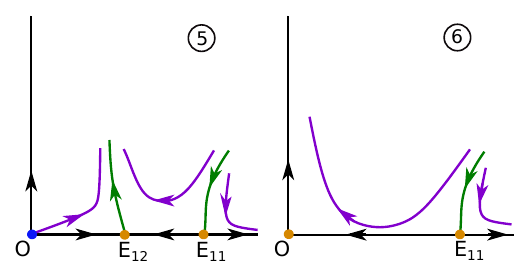} %
\includegraphics[width=0.45\textwidth]{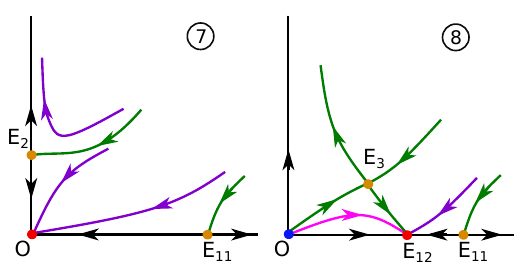} %
\includegraphics[width=0.45\textwidth]{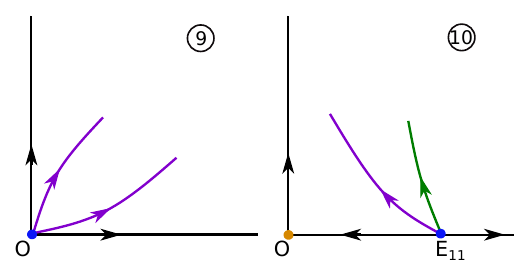} %
\includegraphics[width=0.45\textwidth]{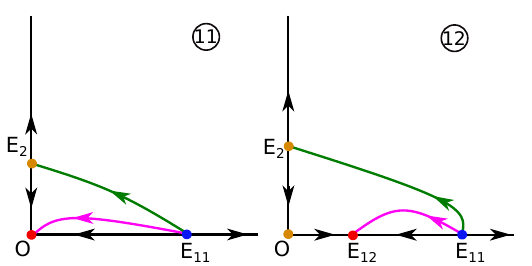} %
\includegraphics[width=0.45\textwidth]{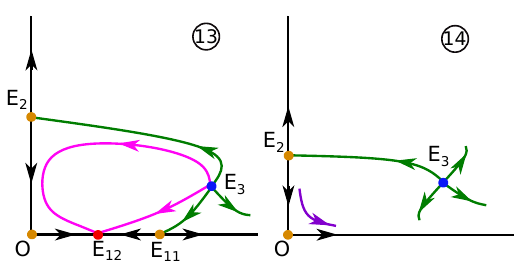} %
\includegraphics[width=0.45\textwidth]{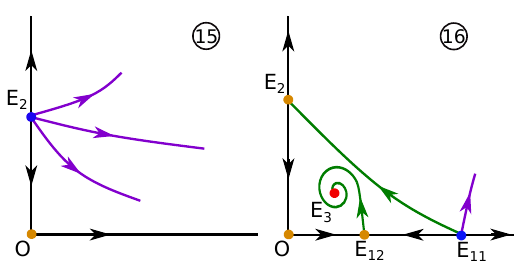} %
\includegraphics[width=0.45\textwidth]{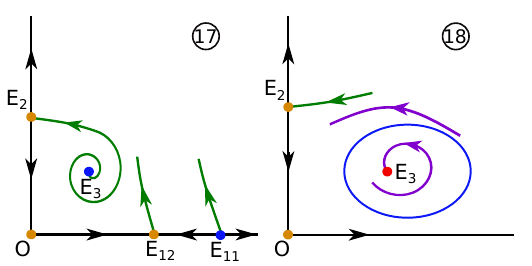} %
\includegraphics[width=0.45\textwidth]{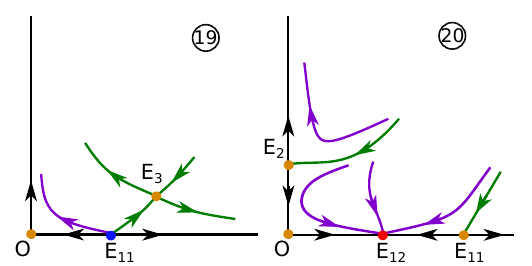} %
\includegraphics[width=0.45\textwidth]{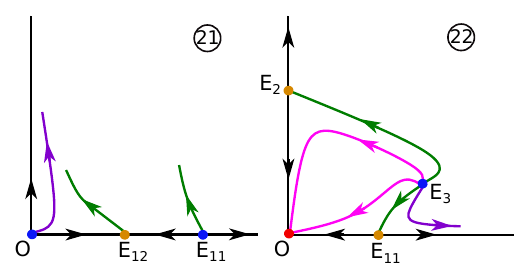}
\caption{The phase portraits corresponding to the diagrams G1-G9}
\label{pp5}
\end{figure}

\begin{figure}[h!]
\centering
\includegraphics[width=0.45\textwidth]{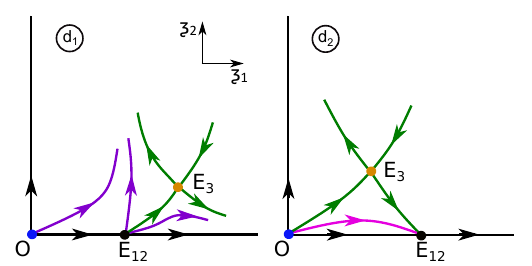} \includegraphics[width=0.45%
\textwidth]{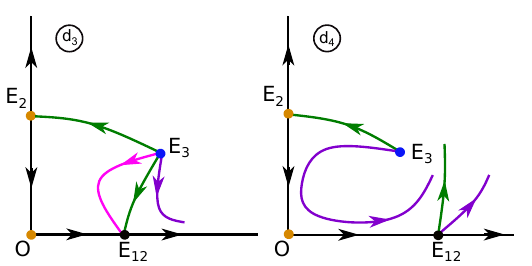} \includegraphics[width=0.45\textwidth]{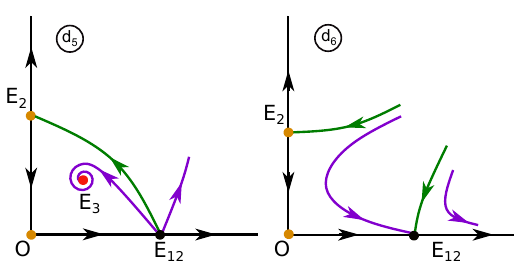} %
\includegraphics[width=0.225\textwidth]{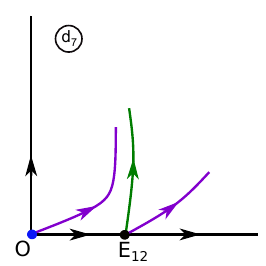}
\caption{The phase portraits for $(\protect\mu_1,\protect\mu_2)\in\Delta$}
\label{pp6}
\end{figure}

\begin{figure}[h!]
\centering
\includegraphics[width=1\textwidth]{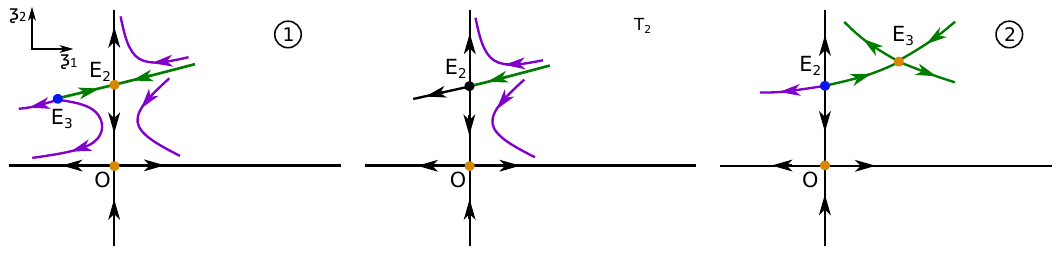}
\caption{The phase portraits on the left and right of $T_2$ in G1}
\label{pp10}
\end{figure}

\begin{figure}[h!]
\centering
\includegraphics[width=1\textwidth]{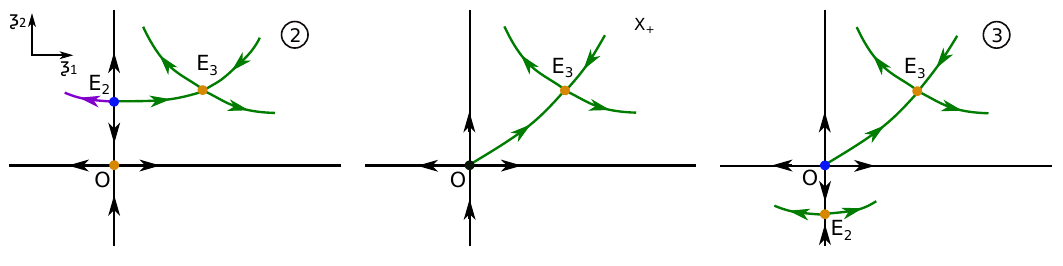}
\caption{The phase portraits on the left and right of $X_+$ in G1}
\label{pp7}
\end{figure}

\begin{figure}[h!]
\centering
\includegraphics[width=1\textwidth]{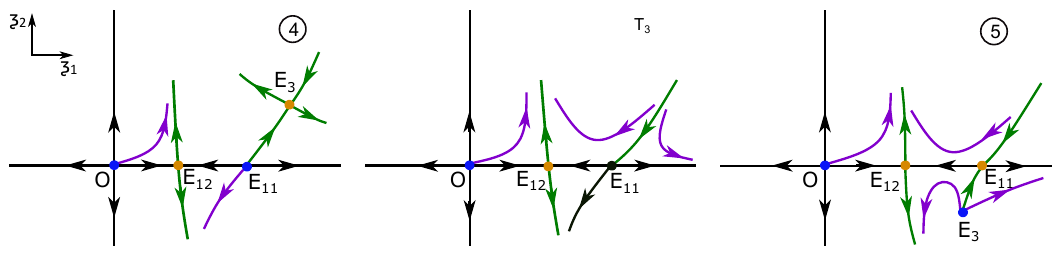}
\caption{The phase portraits on the left and right of $T_3$ in G1}
\label{pp8}
\end{figure}

\begin{figure}[h!]
\centering
\includegraphics[width=1\textwidth]{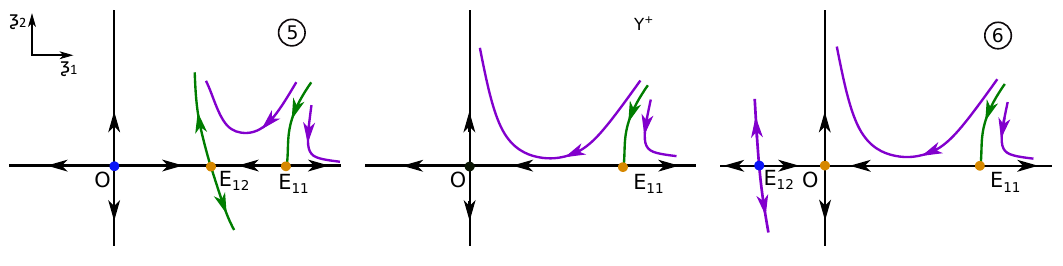}
\caption{The phase portraits on the left and right of $Y_+$ in G1}
\label{pp9}
\end{figure}

\begin{theorem}
Assume $N>0$ and $\sigma _{1}>0.$ Consider the sets 
\begin{equation*}  \label{R10+}
R_{10}^+=\{(\mu_1,\mu_2)\in\mathbb{R}^2\mid
\mu_1>\gamma\sigma_{1}\mu_2^2,\mu_1N<0,\delta\mu_2>0\}\mbox{ and }
R_{10}^-=R_{10}\setminus(R_{10}^+\cup T_3).
\end{equation*}

1) Let $\delta >0.$ If $\theta _{2}>0,$ then $E_{11}$ is a repeller on $%
R_{20}\cup Y_{+}\cup R_{10}^{+},$ a saddle on $R_{10}^{-},$ and unstable on $%
T_{3},$ while $E_{12}$ is a saddle on $R_{20}.$ If $\theta _{2}<0,$ then $%
E_{11}$ is a saddle on $R_{20}\cup Y_{-}\cup R_{10}^{-},$ a repeller on $%
R_{10}^{+},$ and unstable on $T_{3},$ while $E_{12}$ is an attractor on $%
R_{20}.$

2) Let $\delta <0.$ If $\theta _{2}>0,$ then $E_{11}$ is a repeller on $%
R_{20}\cup Y_{+}\cup R_{10}^{-},$ a saddle on $R_{10}^{+},$ and unstable on $%
T_{3},$ while $E_{12}$ is a saddle on $R_{20}.$ If $\theta _{2}<0,$ then $%
E_{11}$ is a saddle on $R_{20}\cup Y_{-}\cup R_{10}^{+},$ a repeller on $%
R_{10}^{-},$ and unstable on $T_{3},$ while $E_{12}$ is an attractor on $%
R_{20}.$
\end{theorem}

\textit{Proof.} By hypothesis we have $N-\delta \theta _{2}>0.$ Using the
notations from Proposition \ref{prop1}, we evaluate $\lambda _{2}^{E_{11}}$
and $\lambda _{2}^{E_{12}}$ on $T_{3}.$ Notice that $E_{3}$ coincides to $%
E_{11}\left( \frac{1}{\delta }\mu _{2},0\right) $ on $T_{3},$ because $%
2N-\delta \theta _{2}>0.$ Furthermore, if $\theta _{2}\mu _{2}>0,$ we have $%
\lambda _{2}^{E_{11}}(0,\mu _{2})=\frac{N-\delta \theta _{2}}{N}\mu
_{2}\left( 1+O(\mu _{2})\right) \neq 0.$

First, if $\delta \mu _{2}>0,$ then $\left. \lambda _{2}^{E_{12}}\right\vert
_{T_{3}}=\frac{2N-\delta \theta _{2}}{N}\mu _{2}\left( 1+O(\mu _{2})\right)
\neq 0.$ Taking into account that $\lambda _{2}^{E_{11}}\lambda
_{2}^{E_{12}}=0$ on $T_{3},$ i.e. $\mu _{1}=\gamma \sigma _{1}\mu _{2}^{2},$
we deduce $\lambda _{2}^{E_{11}}=0$ on $T_{3}.$

1) Let $\delta >0$ and $\theta _{2}>0.$ If $\mu _{2}>0,$ then $\lambda
_{2}^{E_{11}}(0,\mu _{2})>0.$ Because $\lambda _{2}^{E_{11}}=0$ only on $%
T_{3},$ using the above results we obtain that $\lambda _{2}^{E_{11}}>0$ on $%
R_{20}\cup Y_{+}\cup R_{10}^{+}$ and $\lambda _{2}^{E_{11}}<0$ on $%
R_{10}^{-},$ from $\lambda _{2}^{E_{11}}\left( \mu _{1},0\right) <0$ by %
\eqref{la02}.

Moreover, $\lambda _{2}^{E_{12}}>0$ on $R_{20}$ because $\lambda
_{2}^{E_{11}}\lambda _{2}^{E_{12}}>0$ on $R_{20}.$ On the other hand, we
know from Proposition \ref{prop1} that $\lambda _{1}^{E_{11}}>0$ and $%
\lambda _{1}^{E_{12}}<0.$ Therefore, the first conclusion follows.

Now, let $\theta _{2}<0.$ Therefore, on $R_{20}$ \eqref{r10}, $\mu _{2}<0.$
Hence $\lambda _{2}^{E_{11}}(0,\mu _{2})<0.$ As above, $\lambda
_{2}^{E_{11}}<0$ on $R_{20}\cup Y_{-}\cup R_{10}^{-},$ and $\lambda
_{2}^{E_{11}}>0$ on $R_{10}^{+},$ and also $\lambda _{2}^{E_{12}}<0$ on $%
R_{20}.$ Consequently, the second conclusion of assertion 1) is proved.

2) The second part of the theorem follows similarly. $\blacksquare $

\begin{remark}
In Fig. \ref{pp10} we exemplify the behavior of the system \eqref{fn2} when $%
(\mu _{1},\mu _{2})$ crosses $T_{2},$ corresponding to the conditions of G1.
Notice that, the phase portraits from quadrant I ( $\xi _{1}\geq 0,$ $\xi
_{2}\geq 0$) when $(\mu _{1},\mu _{2})\in T_{2}$ and $(\mu _{1},\mu _{2})$
lies in the region "1", i.e. in the region where $E_{3}$ is virtual,
coincide. This occurs because $T_{2}$ is a transcritical bifurcation curve
and the phase portraits are restricted to quadrant I. Similar scenarios take
place when $(\mu _{1},\mu _{2})$ crosses the other transcritical bifurcation
curves, $X_{+}\cup X_{-},$ $T_{3}$ and $Y_{+}\cup Y_{-},$ Figs. \ref{pp7}-%
\ref{pp9}. Restricted to quadrant I, the phase portraits on the bifurcation
curves coincide to the phase portraits corresponding to the regions where
one of the collinding points became virtual after collision.
\end{remark}

The bifurcation diagrams corresponding to the cases $\sigma_1<0$ and $%
\sigma_1>0$ are depicted in Figure \ref{pp3} and Figure \ref{pp4}
respectively.

All possible types of the four equilibrium points arising in the diagrams
are summarized in Table \ref{tab1} and Table \ref{table2}.

\begin{table}[h!]
\centering
\begin{tabular}{l|lllllllllll}
& $1$ & $2$ & $3$ & $4$ & $5$ & $6$ & $7$ & $8$ & $9$ & $10$ & $11$ \\ \hline
$O$ & $s$ & $s$ & $r$ & $r$ & $r$ & $s$ & $a$ & $r$ & $r$ & $s$ & $a$ \\ 
$E_{11}$ & $-$ & $-$ & $-$ & $r$ & $s$ & $s$ & $s$ & $s$ & $-$ & $r$ & $r$
\\ 
$E_{12}$ & $-$ & $-$ & $-$ & $s$ & $s$ & $-$ & $-$ & $a$ & $-$ & $-$ & $-$
\\ 
$E_{2}$ & $s$ & $r$ & $-$ & $-$ & $-$ & $-$ & $s$ & $-$ & $-$ & $-$ & $s$ \\ 
$E_{3}$ & $-$ & $s$ & $s$ & $s$ & $-$ & $-$ & $-$ & $s$ & $-$ & $-$ & $-$%
\end{tabular}%
\caption{The behavior of equilibrium points on different regions from the
bifurcation diagrams G1-G9.}
\label{tab1}
\end{table}

\begin{table}[h]
\centering
\begin{tabular}{l|lllllllllll}
& $12$ & $13$ & $14$ & $15$ & $16$ & $17$ & $18$ & $19$ & $20$ & $21$ & $22$
\\ \hline
$O$ & $s$ & $s$ & $s$ & $s$ & $s$ & $s$ & $s$ & $s$ & $s$ & $r$ & $a$ \\ 
$E_{11}$ & $r$ & $s$ & $-$ & $-$ & $r$ & $r$ & $-$ & $r$ & $s$ & $r$ & $s$
\\ 
$E_{12}$ & $a$ & $a$ & $-$ & $-$ & $s$ & $s$ & $-$ & $-$ & $a$ & $s$ & $-$
\\ 
$E_{2}$ & $s$ & $s$ & $s$ & $r$ & $s$ & $s$ & $s$ & $-$ & $s$ & $-$ & $s$ \\ 
$E_{3}$ & $-$ & $r$ & $r$ & $-$ & $a$ & $r$ & $a$ & $s$ & $-$ & $-$ & $r$%
\end{tabular}%
\caption{Continuation of Table \protect\ref{tab1}.}
\label{table2}
\end{table}

\begin{remark}
We notice that, all bifurcation diagrams G1-G9 of this case do not exist in
the non-degenerate framework. They are new and emerging mainly due to the
existence of the saddle-node bifurcation curves $\Delta _{+}$ and $\Delta
_{-}.$
\end{remark}

\section{Analysis of the system when $\protect\theta \left( 0\right) \neq 0$
and $\protect\delta \left( 0\right) =0$}

Assume $\delta \left( 0\right) =0$ and $\theta =\theta \left( 0\right) \neq
0;$ we still have $\gamma =\gamma \left( 0\right) <0.$ Write in this case $%
\delta \left( \mu \right) =\frac{\partial \delta \left( 0\right) }{\partial
\mu _{1}}\mu _{1}+\frac{\partial \delta \left( 0\right) }{\partial \mu _{2}}%
\mu _{2}+O\left( \left\vert \mu \right\vert ^{2}\right) $ and assume $\frac{%
\partial \delta \left( 0\right) }{\partial \mu _{1}}\overset{not}{=}\delta
_{1}\neq 0.$ The equilibrium points are $O(0,0),$ $E_{1}\left( \frac{1}{%
\theta }\mu _{1}+\allowbreak O\left( \mu _{1}^{2}\right) ,0\right) $ and $%
E_{2}\left( 0,-\mu _{2}+\allowbreak O\left( \mu _{2}^{2}\right) \right) ,$
respectively, 
\begin{equation}
E_{3}\left( \frac{1}{\theta }\left( \mu _{1}-\gamma \mu _{2}\right) ,-\mu
_{2}+\sigma _{2}\mu _{1}^{2}\right)  \label{e3}
\end{equation}%
in its lowest terms, where $\sigma _{2}=\frac{1}{\theta ^{2}}\left( \theta
\delta _{1}-S\right) .$ $E_{3}$ is well-defined and nontrivial for $%
\left\vert \mu \right\vert $ sufficiently small, in the region 
\begin{equation}
Q=\left\{ \left( \mu _{1},\mu _{2}\right) ,\left( \mu _{1}-\gamma \mu
_{2}\right) \theta >0,-\mu _{2}+\sigma _{2}\mu _{1}^{2}>0\right\} .
\label{regq}
\end{equation}%
$E_{3}$ bifurcates from $O$ through two bifurcation curves, namely $T_{2}$
given by \eqref{T2}, $\mu _{2}=\frac{1}{\gamma }\mu _{1}+O\left( \mu
_{1}^{2}\right) , $ $\mu _{1}>0,$ and 
\begin{equation}
T_{4}=\left\{ \left( \mu _{1},\mu _{2}\right) \in 
%TCIMACRO{\U{211d} }%
%BeginExpansion
\mathbb{R}
%EndExpansion
^{2}\mid \mu _{2}=\sigma _{2}\mu _{1}^{2}+O\left( \mu _{1}^{3}\right)
,\theta \mu _{1}>0\right\} .  \label{T4}
\end{equation}%
On $T_{2},$ $E_{3}$ coincides to $E_{2}\left( 0,-\mu _{2}\right) $ while on $%
T_{4}$ to $E_{1}\left( \frac{1}{\theta }\mu _{1},0\right) .$

\begin{remark}
The eigenvalues of $E_{1}\left( \frac{1}{\theta }\mu _{1},0\right) $ are $%
-\mu _{1}$ and $\mu _{2}-\sigma _{2}\mu _{1}^{2},$ respectively of $%
E_{2}\left( 0,-\mu _{2}\right) $ they are $-\mu _{2}$ and $\mu _{1}-\gamma
\mu _{2}.$
\end{remark}

\begin{theorem}
\label{the1} Assume $\theta \delta _{1}\neq 0.$ Then $T_{4}$ is a
transcritical bifurcation curve.
\end{theorem}

\textit{Proof.} When $\left( \mu _{1},\mu _{2}\right) \in T_{4},$ $E_{3}$
coincides to $E_{1}\left( \frac{1}{\theta }\mu _{1},0\right) $ and have the
eigenvalues $0$ and $-\mu _{1}.$ Let $\mu _{2}$ be the bifurcation parameter
while $\mu _{1}\neq 0$ is assumed fixed. For $\mu _{1}\theta >0,$ denote by $%
\xi _{0}=\left( \frac{1}{\theta }\mu _{1},0\right) $ and $\mu _{0}=\left(
\mu _{1},\mu _{2}\right) \in T_{4}.$ Write the system \eqref{fn2} in the
form 
\begin{equation}
\dot{\xi}=f\left( \xi ,\mu \right) ,  \label{ec4}
\end{equation}%
where $\xi =\left( \xi _{1},\xi _{2}\right) ,$ $f=\left( f_{1},f_{2}\right) $
and $\mu =\left( \mu _{1},\mu _{2}\right) .$ On $T_{4},$ the coordinates of $%
E_{3}\left( \xi _{1},\xi _{2}\right) $ satisfy $\xi _{2}=0,$ $\mu
_{2}-\delta \left( \mu \right) \xi _{1}+S\left( \mu \right) \xi _{1}^{2}=0$
and $\mu _{1}-\theta \left( \mu \right) \xi _{1}+N\left( \mu \right) \xi
_{1}^{2}=0.$

Then $A$ and $A^{T}$ have both $\lambda =0$ as an eigenvalue with the
corresponding eigenvectors: $v=\left( 
\begin{array}{c}
\frac{\gamma }{\theta } \\ 
1%
\end{array}%
\right) $ for $A$ and $w=\left( 
\begin{array}{c}
0 \\ 
1%
\end{array}%
\right) $ for $A^{T};$ $A=Df\left( \xi _{0},\mu _{0}\right) $ is the
Jacobian of \eqref{ec4} written in its lowest terms.

In order to write properly $f_{\mu _{2}}=\allowbreak \left( 
\begin{array}{cc}
\frac{\partial f_{1}}{\partial \mu _{2}} & \frac{\partial f_{2}}{\partial
\mu _{2}}%
\end{array}%
\right) ^{T},$ we need $\theta _{2}=\frac{\partial \theta }{\partial \mu _{2}%
}\left( 0\right) ,$ $\delta _{2}=\frac{\partial \delta }{\partial \mu _{2}}%
\left( 0\right) ,$ $M_{2}=\frac{\partial M}{\partial \mu _{2}}\left(
0\right) $ and so on. Then $\frac{\partial f_{1}}{\partial \mu _{2}}=\xi
_{1}\left( -\theta _{2}\xi _{1}+\gamma _{2}\xi _{2}-M_{2}\xi _{1}\xi
_{2}+N_{2}\xi _{1}^{2}\right) $ and $\frac{\partial f_{2}}{\partial \mu _{2}}%
=\xi _{2}\left( S_{2}\xi _{1}^{2}-\delta _{2}\xi _{1}+P_{2}\xi
_{2}^{2}+1\right) ,$ in their lowest terms in $\mu _{1}$ and $\mu _{2}.$

These yield $w^{T}f_{\mu _{2}}\left( \xi _{0},\mu _{0}\right) =0,$ $w^{T}%
\left[ Df_{\mu _{2}}\left( \xi _{0},\mu _{0}\right) \left( v\right) \right]
=1+O\left( \mu _{1}\right) \neq 0$ and $w^{T}\left[ D^{2}f\left( \xi
_{0},\mu _{0}\right) \left( v,v\right) \right] =2+O\left( \mu _{1}\right)
\neq 0,$ for $\left\vert \mu _{1}\right\vert $ small. Thus, $T_{4}$ is a
transcritical bifurcation curve. $\blacksquare $

\begin{remark}
$T_{2},$ $Y_{+},$ $Y_{-},$ $X_{+}$ and $X_{-}$ are transcritical bifurcation
curves. The proof is similar to Theorem \ref{the1}.
\end{remark}

\begin{theorem}
Assume $\theta \delta _{1}\neq 0$ and $\left( \mu _{1},\mu _{2}\right) \in
Q. $ Then, $E_{3}$ is a saddle when $\theta >0,$ respectively, a repeller
when $\theta <0.$ A Hopf bifurcation cannot occur at $E_{3}.$
\end{theorem}

Proof. The characteristic polynomial at $E_{3}$ is $P\left( \lambda \right)
=\lambda ^{2}-2p\lambda +L,$ where $p$ and $L$ are given by \eqref{pin} and %
\eqref{q1}; $L=-\xi _{1}\xi _{2}\left( \theta +O\left( \left\vert \mu
\right\vert \right) \right) .$ The eigenvalues $\lambda _{1,2}=p\pm \sqrt{q}$
at $E_{3}$ satisfy $sign\left( \lambda _{1}\lambda _{2}\right) =-sign\left(
\theta \right) .$ Thus, $E_{3}$ is a saddle if $\theta >0.$

If $\theta <0,$ then $p=0$ is a curve given by 
\begin{equation}
H_{1}=\left\{ \left( \mu _{1},\mu _{2}\right) ,\mu _{2}=\frac{1}{\gamma -1}%
\mu _{1}+O\left( \mu _{1}^{2}\right) \right\} .  \label{h2}
\end{equation}%
In its lowest terms, $p$ reads $p=\frac{1}{2}\left( \gamma -1\right)
\allowbreak \mu _{2}-\frac{1}{2}\mu _{1}$ and is given by $p=-\frac{1}{2}%
\theta \xi _{1}+\frac{1}{2}\xi _{2}>0,$ whenever $E_{3}$ exists. It follows
that $H_{1}\nsubseteq Q$ and $p>0$ on $Q,$ thus, $E_{3}$ is a repeller and a
Hopf bifurcation cannot occur at $E_{3}.$ Notice that $\gamma \neq 1$
because $\gamma <0.$ $\blacksquare $

\begin{figure}[tbp]
\centering
\includegraphics[width=0.8\textwidth]{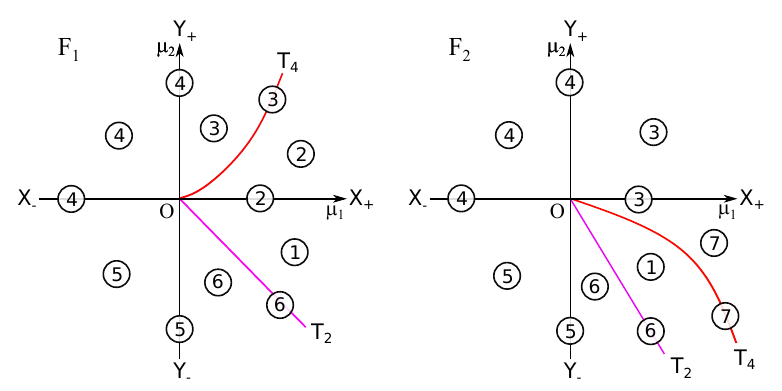} \includegraphics[width=0.8%
\textwidth]{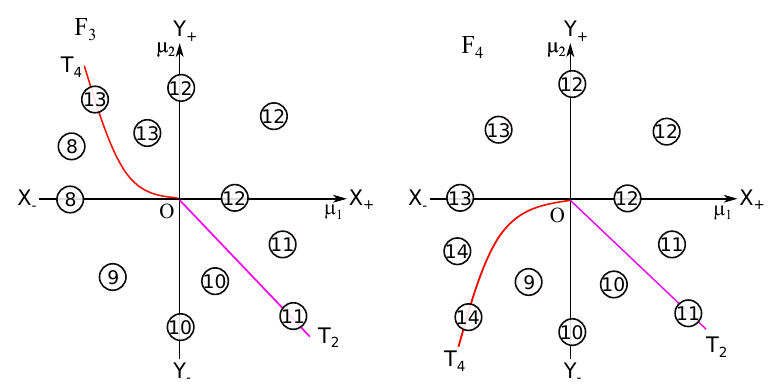}
\caption{Bifurcation diagrams corresponding to $\protect\theta \left( 0\right) \neq 0$
and $\protect\delta \left( 0\right) =0.$ }
\label{dia1}
\end{figure}

\begin{figure}[tbp]
\centering
\includegraphics[width=0.45\textwidth]{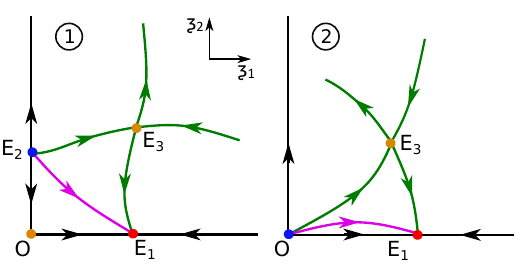} %
\includegraphics[width=0.45\textwidth]{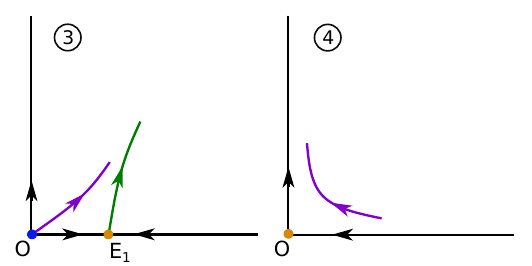} %
\includegraphics[width=0.45\textwidth]{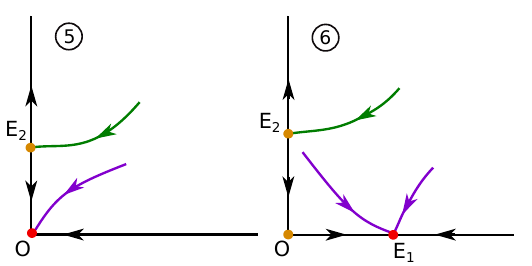} %
\includegraphics[width=0.45\textwidth]{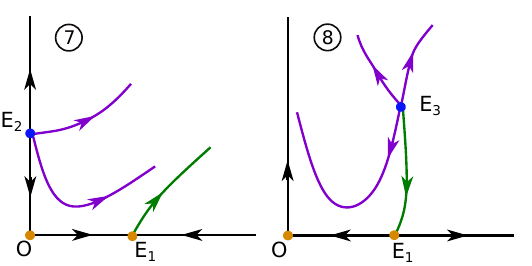} %
\includegraphics[width=0.45\textwidth]{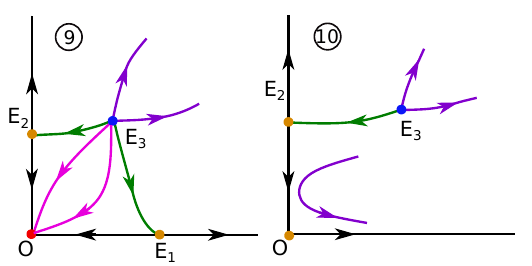} %
\includegraphics[width=0.45\textwidth]{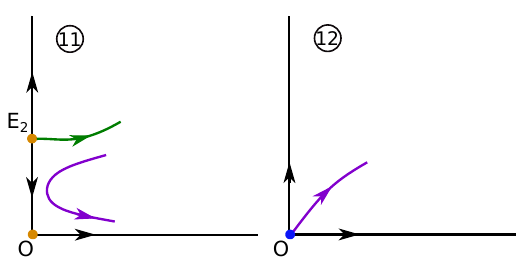} %
\includegraphics[width=0.45\textwidth]{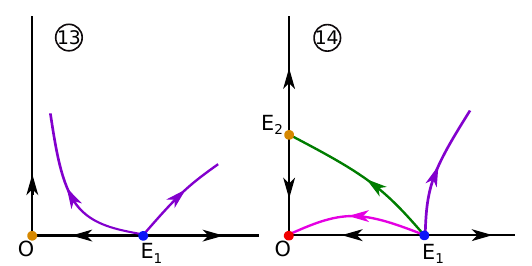}
\caption{The phase portraits corresponding to the diagrams $F_1-F_4.$}
\label{dia2}
\end{figure}

\begin{remark}
The bifurcation diagrams F1-F4 corresponding to this degeneracy present
similarities with four diagrams from the non-degenerate case, namely with
the diagrams VII-X reported in \cite{GT1}. However, two bifurcation curves
from these diagrams are different in the two cases, namely $T_{4}$ from the
degenerate case is a parabola-like curve, while $T_{1}$ is a line in the
non-degenerate framework. Moreover, the equilibrium points have different
expressions in the two cases.
\end{remark}

\section{Conclusions}

We approached in this work two degenerate cases, namely $\theta \left(
0\right) =0$ and $\delta \left( 0\right) \neq 0,$ respectively, $\theta
\left( 0\right) \neq 0$ and $\delta \left( 0\right) =0.$ Another degenerate
case, $\theta \left( 0\right) =0$ and $\delta \left( 0\right) =0,$ is more
involved due to the fact that the Implicit Function Theorem cannot be
applied anymore for finding the equilibrium $E_{3}.$ The behavior of the
system in this case is an open problem.

\section{Acknowledgments}

This research was supported by Horizon2020-2017-RISE-777911 project. We
thank to Prof. Jaume Llibre for his useful suggestions related to Kolmogorov
systems. We are also grateful to the anonymous reviewers for their comments which improved the paper.

\end{document}